\documentclass{amsart}
\usepackage{amssymb,amscd,enumerate,xspace,amsthm}
\usepackage[all,ps,cmtip]{xy}

\theoremstyle{plain}
\newtheorem{theoreme}[equation]{Theorem}
\newtheorem{proposition}[equation]{Proposition}
\newtheorem{corollaire}[equation]{Corollary}
\newtheorem{lemme}[equation]{Lemma}
\newtheorem*{assertion*}{Assertion}

\theoremstyle{remark}
\newtheorem{Remarques}[equation]{Remarks}
\newtheorem{remarque}[equation]{Remark}
\newtheorem*{remarque*}{Remark}
\newtheorem*{Remarques*}{Remarks}
\newtheorem{definition}[equation]{Definition}

\numberwithin{equation}{section}

\let\dpl\displaystyle

\let\wt\widetilde
\let\ov\overline

\newcommand{\bbullet}{{\scriptscriptstyle\bullet}}

\newcommand{\Implique}{\Longrightarrow}

\newcommand{\qbox}[1]{\quad\hbox{#1}}
\newcommand{\qqbox}[1]{\quad\hbox{#1}\quad}

\newcommand{\ooplus}{\mathop\oplus\limits}
\newcommand{\ootimes}{\mathop\otimes\limits}

\newcommand{\module}[1]{\left\vert#1\right\vert}

\DeclareMathOperator{\coker}{Coker}
\newcommand{\defin}{\stackrel{\mathrm{def}}{=}}

\DeclareMathOperator{\gr}{gr}
\DeclareMathOperator{\id}{Id}
\DeclareMathOperator{\im}{Im}
\newcommand{\isom}{\stackrel{\sim}{\longrightarrow}}

\renewcommand{\ker}{\mathop{\rm Ker}\nolimits}

\newcommand{\reel}{\mathop{\mathrm{Re}}\nolimits}
\newcommand{\res}{\mathop{\mathrm{Res}}\nolimits}

\newcommand{\cHom}{\mathop{\cH om}\nolimits}
\newcommand{\rHom}{\mathop{\bR\cH om}\nolimits}
\newcommand{\ext}{\mathop{\cE xt}\nolimits}

\DeclareMathOperator{\supp}{Supp}

\newcommand{\loccit}{\emph{loc\ptbl cit}}
\newcommand{\eg}{{\it e.g}}
\newcommand{\cf}{\emph{cf}}
\newcommand{\etc}{{\it etc}}
\newcommand{\ie}{{\it i.e}}
\newcommand{\resp}{{\it resp}}

\newcommand{\T}{\S\kern .15em }
\newcommand{\ptbl}{.\kern .15em }

\newcommand{\lefpar}{\left(}
\newcommand{\rigpar}{\right)}
\newcommand{\lefcro}{\left[}
\newcommand{\rigcro}{\right]}

\newdimen\lengtharrow
 \newbox\exponantbox \newbox\indicebox
\def\dimmax#1#2{\ifdim#1<#2 #2\else #1\fi}

\def\arrowr#1#2%
    {\setbox\exponantbox=\hbox{$#2$}
      \setbox\indicebox=\hbox{$#1$}
       \lengtharrow=\dimmax{\wd\indicebox}{\wd\exponantbox}
        \ifdim \lengtharrow=0pt
                                                                                                                        \mathrel{\smash{\mathop{\hbox to 10mm{\rightarrowfill}}}}
         \else \ifdim\lengtharrow<6mm
                \lengtharrow=10mm
                 \else \advance\lengtharrow by 4mm
                  \fi
                                        \mathrel{\smash{\mathop{\hbox to\lengtharrow{\rightarrowfill}}\limits%
                         ^{\styleflechehoriz#2}_{\styleflechehoriz#1}}}
                                                                        \fi}

\def\arrowl#1#2%
    {\setbox\exponantbox=\hbox{$#2$}
      \setbox\indicebox=\hbox{$#1$}
       \lengtharrow=\dimmax{\wd\indicebox}{\wd\exponantbox}
        \ifdim \lengtharrow=0pt
                                                                                                                        \mathrel{\smash{\mathop{\hbox to 10mm{\leftarrowfill}}}}
         \else \ifdim\lengtharrow<6mm
                \lengtharrow=10mm
                 \else \advance\lengtharrow by 4mm
                  \fi
                                        \mathrel{\smash{\mathop{\hbox to\lengtharrow{\leftarrowfill}}\limits%
                         ^{\styleflechehoriz#2}_{\styleflechehoriz#1}}}
                                                                        \fi}

\def\hookarrowr#1#2{\lhook\joinrel\arrowr{#1}{#2}}

\def\dasharrowr#1#2{\mapstochar\arrowr{#1}{#2}}

\def\MRE#1{\arrowr{}{#1}}

\def\HMRE#1{\hookarrowr{}{#1}}

\def\DMRE#1{\dasharrowr{}{#1}}

\let\styleflechehoriz=\textstyle

\def\CC{\mathbf{C}}

\def\NN{\mathbf{N}}

\def\RR{\mathbf{R}}

\def\ZZ{\mathbf{Z}}

\def\cA{\mathcal{A}}

\def\cC{\mathcal{C}}
\def\cD{\mathcal{D}}
\def\cE{\mathcal{E}}
\def\cF{\mathcal{F}}

\def\cH{\mathcal{H}}

\def\cJ{\mathcal{J}}

\def\cM{\mathcal{M}}
\def\cN{\mathcal{N}}
\def\cO{\mathcal{O}}

\def\cR{\mathcal{R}}

\def\cT{\mathcal{T}}

\def\bL{\boldsymbol{L}}

\def\bR{\boldsymbol{R}}

\DeclareMathAlphabet{\mathcalmaigre}{U}{eus}{m}{n}

\def\ccS{\mathcalmaigre{S}}

\DeclareMathAlphabet{\mathcalgras}{U}{eus}{b}{n}

\def\eg{\emph{e.g}}
\newcommand{\cTor}{\mathop{\cT\!or}\nolimits}
\newcommand{\cbbullet}{\raisebox{1pt}{$\bbullet$}}
\DeclareMathOperator{\can}{can}

\DeclareMathOperator{\Var}{Var}
\DeclareMathOperator{\coloc}{coloc}
\DeclareMathOperator{\loc}{loc}

\DeclareMathOperator{\DR}{DR}
\DeclareMathOperator{\Mod}{Mod}

\DeclareMathOperator{\Db}{\mathfrak{Db}}
\DeclareMathOperator{\gC}{\mathfrak{C}}

\DeclareMathOperator{\RHDb}{RH\Db}
\DeclareMathOperator{\CRHDb}{\cC^\infty\! RH\Db}

\def\Dbmo{\Db^{{\rm mod}\,0}}

\def\RHDbmo{\RHDb^{{\rm mod}\,0}}
\def\Dbmz{\Db^{{\rm mod}\,Z}}

\def\RHDbmz{\RHDb^{{\rm mod}\,Z}}

\def\CRHDbmo{\CRHDb^{{\rm mod}\,0}}

\def\psim{\psi^{\rm mod}}
\def\phim{\phi^{\rm mod}}
\def\Cmz{C^{{\rm mod}\,Z}}

\let\dag\dagger

\def\loc{{\rm loc}}

\let\la\llangle
\let\ra\rrangle

\begin{document}

\title{Vanishing cycles and Hermitian duality}

\author{Claude Sabbah}
\address{UMR 7640 du CNRS\\
Centre de Math{\'e}matiques\\
{\'E}cole polytechnique\\
F--91128 Palaiseau cedex\\
France}
\email{sabbah@math.polytechnique.fr}
\urladdr{http://math.polytechnique.fr/cmat/sabbah/sabbah.html}
\thanks{This work is partly supported by INTAS program 97-1644}

\begin{abstract}
We show the compatibility between the moderate or nearby cycle functor for regular holonomic $\mathcal{D}$-modules, as defined by Beilinson, Kashiwara and Malgrange, and the Hermitian duality functor, as defined by Kashiwara.
\end{abstract}
\maketitle

\section*{Introduction}
The Hermitian dual of a $\cD$-module was introduced by M.~Kashiwara in \cite{Kashiwara86}, who showed that the Hermitian dual of a regular holonomic $\cD$-module is also regular holonomic (hence coherent). In this paper we show a compatibility result between this functor and the nearby or vanishing cycle functor relative to a holomorphic function for such modules. The latter may be defined using the $V$-filtration (introduced by Beilinson, Kashiwara and Malgrange). 

Moreover we make the link with asymptotic expansions of integrals along fibres of the function. This gives a generalization of previous work of D.~Barlet on Hermitian duality for the local Gauss-Manin system of an analytic function. In particular this gives a simpler approach to the ``tangling phenomenon'' described by D.~Barlet in \cite{Barlet91}.

\subsubsection*{Acknowledgements}
This work arose from many discussions with D.~Barlet, whom I thank.

\section{Hermitian duality}

\subsection{Notation}
Let $(X,\cO_X)$ be a complex analytic manifold of dimension $n$, $(X_\RR,\cA_{X_\RR})$ be the underlying real analytic manifold and let $(\ov X,\cO_{\ov X}=\ov\cO_X)$ be the complex conjugate manifold. Denote by $\cD_X$ (\resp.\ $\cD_{\ov X}$) the sheaf of holomorphic linear differential operators on $X$ (\resp.\ $\ov X$).

Denote by $\ov{\phantom{X}}:f\mapsto \ov f$ the $\RR$-isomorphism $\cO_X\rightarrow \cO_{\ov X}$ and $\cD_X\rightarrow \cD_{\ov X}$. It induces a trivial conjugation functor, sending $\cD_X$-modules to $\cD_{\ov X}$-modules; if $\cM$ is a $\cD_X$-module, we denote by $\ov\cM$ the $\RR$-vector space $\cM$ equipped with the action of $\cD_{\ov X}$ defined as follows: denote by $\ov m$ the local section $m$ of $\cM$ viewed as a local section of $\ov M$; then $\ov P\cdot \ov m= Pm$.

Let $\Db_{X_\RR}$ (also denoted by $\Db_X$ for short) be the sheaf of distributions on $X_\RR$. It acts on the sheaf $C^\infty$-forms $\varphi$ with compact support of maximal degree, which is a right $\cD_X$ and $\cD_{\ov X}$-module. Then $\Db_X$ is a left $\cD_X$ and $\cD_{\ov X}$-module by the formula $(P\ov Q\mu)(\varphi)=\mu(\varphi\cdot P\ov Q)$. The sheaf $\gC_{X_\RR}=\Db_{X}^{(n,n)}$ of currents of maximal degree is a right $\cD_X$ and $\cD_{\ov X}$-module obtained from $\Db_X$ by ``going from left to right''.

It will be convenient in the following to denote by $\cO_{X,\ov X}$ (\resp.\ $\cD_{X,\ov X}$) the sheaf $\cO_X\otimes_\CC\cO_{\ov X}$ (\resp.\ $\cD_X\otimes_\CC\cD_{\ov X}$) and to view $\Db_X$ (\resp.\ $\Db_{X}^{(n,n)}$) as a left (\resp.\ right) $\cD_{X,\ov X}$-module.

Let $Z$ be a reduced divisor in $X$ and $\cO_X[*Z]$ the sheaf of meromorphic functions on $X$ with poles along $Z$. There is an exact sequence of left $\cD_{X,\ov X}$-modules
\[
0\longrightarrow \Db_{X,Z}\longrightarrow \Db_X\longrightarrow \Dbmz_X\longrightarrow 0
\]
where $\Db_{X,Z}$ denotes the sheaf of distributions supported on $Z_\RR$ and (see \eg. \cite[Chap\ptbl VII]{Malgrange66})
\[
\Dbmz_X=\cO_X[*Z]\ootimes_{\cO_X}\Db_X={\rm image}\lefcro \Db_X\rightarrow j_*\Db_{X-Z}\rigcro
\]
denotes the subsheaf of $j_*\Db_{X-Z}$ (where $j:X-Z\hookrightarrow X$ denotes the open inclusion) of distributions on $X-Z$ with {\em moderate growth} along $Z$.

\subsection{The Hermitian duality functor \protect\cite{Kashiwara86}}
Denote by $C_X$ the Hermitian duality functor\footnote{It is called improperly the ``conjugation functor'' in \cite{Bjork93}.}. Recall that $C_X$ is a contravariant functor from the derived category $D^-(\cD_X)$ to the category $D^+(\cD_{\ov X})$ defined as
\begin{eqnarray*}
C_X(\cM^\bbullet)&=&\rHom_{\cD_X}(\cM^\bbullet,\Db_X).
\end{eqnarray*}
It restricts as a functor from the full subcategory $D^b_{hr}(\cD_X)$ of bounded complexes with regular holonomic cohomology to $D^b_{hr}(\cD_{\ov X})$ and is equal to the functor $\cHom_{\cD_X}(\bbullet,\Db_X)$ on the category of regular holonomic $\cD_X$-modules (see \cite{Kashiwara86}, see also \cite[Chap\ptbl VII]{Bjork93}), defining there an anti-equivalence of categories between $\Mod_{hr}(\cD_X)$ and $\Mod_{hr}(\cD_{\ov X})$, and between $D^b_{hr}(\cD_X)$ and $D^b_{hr}(\cD_{\ov X})$, $C_{\ov X}$ being a quasi-inverse functor. On $D^b_{hr}(\cD_X)$ we have
\begin{eqnarray*}
\cH^kC_X\cM^\bbullet&=&C_X\cH^k\cM^\bbullet.
\end{eqnarray*}

Last, recall (see \cite{Kashiwara86}) that the \emph{conjugate} of a regular holonomic right $\cD_X$-module $\cM$ is the right $\cD_X$-module defined as
\begin{eqnarray*}
\cM^c&=&\cTor_n^{\cD_{\ov X}}(\ov\cM,\Db_{X}^{(n,0)})
\end{eqnarray*}
and satisfies $\DR\cM^c=\ov{\DR\cM}$. The conjugate of a left module is then obtained in the usual way.

For $Z$ as above, we will denote by $\Cmz_X$ the functor defined as
\[
\begin{array}{ccccc}
\Cmz_X(\cM)&=&C_X(\cM)[*\ov Z]&=&\cHom_{\cD_X}(\cM,\Dbmz_X).
\end{array}
\]
We will call $\ov{C_X(\cM)}$ the \emph{Hermitian dual} of $\cM$.

We may also define the Hermitian dual of a right regular holonomic $\cD$-module by using the sheaf $\Db_{X}^{(n,n)}$ of currents instead of the sheaf $\Db_X$ of distributions.

\begin{remarque}[Extension to the holonomic case]
Kashiwara conjecured (see \cite[Rem\ptbl3.5]{Kashiwara86}) that the previous results remain true for holonomic modules. This is proved in \cite{Bibi97} when the support of $\cM$ has dimension $1$ and in some cases when it has dimension $2$. If this conjecture is true, Theorem \ref{th:compatible} also applies to holonomic modules. It would then be interesting to extend Theorem \ref{th:bm} below to the nonregular holonomic case in order to get a holonomic analogue of Theorem \ref{th:cxa}.
\end{remarque}

\subsection{Sesquilinear forms on $\cD_X$-modules}
Let $\cM',\cM''$ be two left $\cD_X$-modules. A sesquilinear form will be a $\cD_{X,\ov X}$-linear morphism
\begin{eqnarray*}
S:\cM'\ootimes_\CC \ov{\cM''}&\longrightarrow &\Db_X.
\end{eqnarray*}
The datum of $S$ is equivalent to the datum of a $\cD_{\ov X}$-linear morphism
\begin{eqnarray*}
L_S:\ov{\cM''}&\longrightarrow& \cHom_{\cD_X}(\cM'\Db_X).
\end{eqnarray*}
We say that $S$ is nondegenerate if this morphism is an isomorphism.

When $\cM'=\cM''=\cM$ is regular holonomic, this is equivalent to saying that $L_S:\ov \cM\rightarrow C_X\cM$ is injective (or surjective), because $\ov\cM$ and $C_X\cM$ have the same characteristic variety (as their de~Rham complexes are Verdier dual one to each other). We say that $S$ is $\pm$-Hermitian if $\ov{C_X(L_S)}=\pm L_S$, in other words if $S(m,\ov \mu)=\pm S(\mu,\ov m)$ in $\Db_X$.

\subsection{Direct and inverse image by a closed immersion of codimension one and Hermitian duality}
Let $Z$ be a reduced divisor in $X$ and $i:Z\hookrightarrow X$ (\resp.\ $\ov\imath:\ov Z\hookrightarrow \ov X$) denote the inclusion. Let $j_+j^+$ be the localization functor along $Z$ and denote by $j_\dag j^+$ its adjoint by duality, \ie.\ $j_\dag j^+=Dj_+j^+D$, where $D$ denotes the duality functor on holonomic $\cD_X$-modules given by $D\cM=\cHom_{\cD_X}\lefpar \Omega_X^n,\ext_{\cD_X}^{n}(\cM,\cD_X)\rigpar$, with $n=\dim X$.

We also consider the two functors $i_+i^+$ and $i_+i^\dag$. Recall that, for a holonomic $\cD_X$-module $\cM$, we have the following two dual exact sequences
\[
\xymatrix @R=5mm {
0\ar[r]&\cH^{-1}(i_+i^+\cM)\ar[r]&\cM\ar[r]^-{\loc}& j_+j^+\cM\ar[r]&\cH^0(i_+i^+\cM)\ar[r]&0\\
0\ar[r]&\cH^0(i_+i^\dag\cM)\ar[r]&j_\dag j^+\cM\ar[r]_-{\coloc}&\cM\ar[r]& \cH^1(i_+i^\dag\cM)\ar[r]&0
}
\]

\begin{proposition}\label{prop:imdir}
There is a natural isomorphism of contravariant functors from $\Mod_{hr}(\cD_X)$ to $\Mod_{hr}(\cD_{\ov X})$
\begin{eqnarray*}
\Cmz_X&\simeq&C_X\circ j_\dag j^+
\end{eqnarray*}
under which $C_X(\coloc_\cM)$ corresponds to $\loc_{C_X\cM}$.
\end{proposition}

\begin{proof}
The first part is proved in \cite[Proposition II.3.2.2]{Bibi97}. We now want to prove that the following diagram
\[
\xymatrix@C=2cm@R=1cm{
C_X(\cM)\ar[r]^-{C_X(\coloc_\cM)}\ar[rrd]_-{\loc_{C_X\cM}}& C_X(j_\dag j^+\cM)\ar[r]^-{\sim}&\Cmz_X(j_\dag j^+\cM)\\
&&\Cmz_X(\cM)\ar[u]^-{\wr}_-{\Cmz_X(\coloc_\cM)}
}
\]
commutes. Remark first that it clearly commutes on $X-Z$. Put $\cN=\ov{C_X(\cM)}$. The upper part of the diagram gives a morphism $\varphi:\cN\rightarrow \cN[*Z]$ which induces the identity on $X-Z$. It thus factorizes uniquely through $\loc:\cN\rightarrow \cN[*Z]$ to give a morphism $\psi:\cN[*Z]\rightarrow \cN[*Z]$ equal to $\id$ on $X-Z$. It follows that $\psi=\id$ (indeed, $\psi$ is injective because $\cN[*Z]$ has no torsion supported on $Z$, and therefore is onto as $\cN[*Z]$ is holonomic).
\end{proof}

\begin{corollaire}\label{cor:imdir}
The nondegenerate pairing $$j_\dag j^+\cM\ootimes_\CC \Cmz_X(\cM)\longrightarrow \Db_X$$ induces a nondgenerate pairing $$\cH^{-k}(\ov\imath_+\ov\imath^+C_X(\cM))\ootimes_\CC\cH^k(i_+i^\dag\cM)\longrightarrow \Db_X$$ and hence an isomorphism $$\cH^{-k}(\ov\imath_+\ov\imath^+C_X(\cM))\isom C_X\cH^k(i_+i^\dag\cM)$$ for $k=0,1$.\hfill\qed
\end{corollaire}

\begin{corollaire}\label{cor:iminv}
Assume that $Z$ is smooth. Then there is a natural isomorphism of functors ($k=0,1$)
\[
C_Z\circ \cH^k(i^\dag)\simeq\cH^{-k}(\ov\imath^+)\circ C_X.
\]
\end{corollaire}

\begin{proof}
Remark first that there is a natural isomorphism of functors
\[
C_X\circ i_+\simeq\ov\imath_+\circ C_Z.
\]
Indeed, denoting by $i_{++}$ the direct image of $\cD_{Z,\ov Z}$-modules, recall that one has $\Db_{X,Z}=i_{++}\Db_Z$: indeed, put $\cD_{X\leftarrow Z}\otimes_\CC\cD_{\ov X\leftarrow \ov Z}=\cD_{Z\leftarrow X,\ov Z\leftarrow\ov X}$ and consider the natural morphism of right $\cD_{X,\ov X}$-modules
\[
\gC_Z\ootimes_{\cD_{Z,\ov Z}}\lefpar \cD_{Z\leftarrow X,\ov Z\leftarrow\ov X}\rigpar \longrightarrow \gC_X
\]
such that, for any local section $\mu$ of $\gC_Z$, the image of $\mu\otimes 1$ evaluated on any function $\varphi\in\cC^\infty_c(X)$ is equal to $\mu(\varphi_{|Z})$; this morphism is an isomorphism, as can be seen from a local computation; going from right to left, one gets the assertion.

It follows that
\begin{align*}
\cHom_{\cD_X}(i_+\cM,\Db_X)&=\rHom_{\cD_X}(i_+\cM,\Db_{X,Z})\\
&=\bR i_*\rHom_{i^{-1}\cD_X}\Bigl(\cD_{X\leftarrow Z}\ootimes_{\cD_Z}^{\bL}\cM,\cD_{X\leftarrow Z,\ov X\leftarrow \ov Z}\ootimes_{\cD_{Z,\ov Z}}^{\bL}\Db_Z\Bigr)\\
&=\bR i_*\cD_{\ov X\leftarrow \ov Z}\ootimes_{\cD_{\ov Z}}^{\bL} \rHom_{i^{-1}\cD_X}\Bigl(\cD_{X\leftarrow Z} \ootimes_{\cD_Z}^{\bL}\cM,\cD_{X\leftarrow Z}\ootimes_{\cD_{Z}}^{\bL}\Db_Z\Bigr)\\
&=\bR i_*\cD_{\ov X\leftarrow \ov Z}\ootimes_{\cD_{\ov Z}}^{\bL} \rHom_{\cD_Z}(\cM,\Db_Z)\quad\text{(Kashiwara's equivalence)}\\
&=\ov\imath_+ C_Z\cM.
\end{align*}
As $i_+$ and $\ov\imath_+$ are exact functors, we obtain from Corollary \ref{cor:imdir} an isomorphism
\[
\ov\imath_+C_Z\circ \cH^k(i^\dag)\simeq\ov\imath_+\cH^{-k}(\ov\imath^+)\circ C_X,
\]
and thus the result, as $\ov\imath_+$ is an equivalence.
\end{proof}

\section{Regular holonomic distributions}
\subsection{Regular holonomic distributions \cite{Kashiwara86}, \cite[chap\ptbl VII]{Bjork93}}
Let $\Omega$ be an open set in $X$. A distribution $u\in\Db(\Omega)$ is {\em regular holonomic} if the sub-$\cD_\Omega$-module $\cD_\Omega\cdot u$ (or equivalently the sub-$\cD_{\ov\Omega}$-module $\cD_{\ov\Omega}\cdot u$, \cf.\ \cite[Proposition 4]{Kashiwara86}, \cite[Proposition 7.4.2]{Bjork93}) of $\Db_\Omega$ is regular holonomic. The notion is local, \ie.\ there exists a sheaf $\RHDb_X$ such that the set of regular holonomic distributions on $\Omega$ is $\Gamma(\Omega,\RHDb_X)$.

Notice that $\RHDb_X$ is a left $\cD_X$ and $\cD_{\ov X}$-module. It will be convenient to consider the subsheaf $\cC^\infty_X\cdot\RHDb_X$ of $\Db_X$ whose local sections are finite combinations of regular holonomic distributions with $C^\infty$ coefficients.

Analogous results hold for $\RHDbmz_X$. Notice that we have
\[
\RHDbmz_X=\cO_X[*Z]\ootimes_{\cO_X}\RHDb_X={\rm image}\lefcro \RHDb_X\rightarrow \Db_{X-Z}\rigcro.
\]

The following is a slight generalization of \cite{Barlet82} and \cite[Theorem 11]{B-M87}.

\begin{theoreme}\label{th:bm}
Let $X=Z\times \CC$ have dimension $n+1$ and let $u$ be a regular holonomic distribution on the open set $\Omega\times D$ of $X$. Let $\varphi\in\cD^{n,n}(\Omega)$ be a $C^\infty$ $(n,n)$-form with compact support. Then $\langle u,\varphi\rangle\in\Db(D)$ is in $\Gamma(D,\CRHDb_\CC)$.
\end{theoreme}

\begin{proof}
According to Remark \ref{rem:dvptasympt} and Proposition \ref{prop:caractRH} below, it is identical to the one of \loccit., using the existence of a good Bernstein relation (\cf.\ for instance \cite[Theorem 8.8.16]{Bjork93} and the references given there) in order to prove the existence of a good operator as in \cite[Proposition 8]{B-M87}.
\end{proof}

\subsection{Regular holonomic distributions in dimension $1$}\label{subsec:RH1}
Assume now that $X=\CC$ and $Z=\{0\}$. Let $t$ be a coordinate on $\CC$. For $\alpha\in\CC$ such that $-1\leq \reel\alpha<0$ and $p\in\ZZ$, put
\begin{eqnarray*}
u_{\alpha,p}&=&
\begin{cases}
\module{t}^{-2(\alpha+1)}\dfrac{(\log\module{t}^2)^p}{p!}\in L^1_{\rm loc}(\CC)& \text{if }p\geq 0\\
0&\text{if }p<0.
\end{cases}
\end{eqnarray*}
then it is easy to show that the family $u_{\alpha,p}$ satisfies
\begin{equation}\label{eq:uap}
(\partial_tt+\alpha)u_{\alpha,p}=(\partial_{\ov t}\ov t+\alpha)u_{\alpha,p}=u_{\alpha,p-1}.
\end{equation}
This implies in particular that $u_{\alpha,p}\in\RHDb_{\CC,0}$.

Let $\RHDbmo_{\CC,0}$ denote the image of $\RHDb_{\CC,0}$ in $\Dbmo_{\CC,0}=\Db_{\CC,0}[1/t]=\Db_{\CC,0}[1/\ov t]$. We then have
\[
\RHDbmo_{\CC,0}=\RHDb_{\CC,0}[1/t]=\RHDb_{\CC,0}[1/\ov t].
\]
It is known (see \cite[chap\ptbl VII, \T7]{Bjork93}) that
\begin{eqnarray}\label{eq:RHDbmo}
\RHDbmo_{\CC,0}&=&\sum_{-1\leq \reel\alpha<0}\sum_{p}\CC\{t\}[t^{-1}]\cdot\CC\{\ov t\}[\ov t^{-1}]\cdot \wt u_{\alpha,p}
\end{eqnarray}
where $\wt u$ denotes the image of the distribution germ $u$ in $\Dbmo_{\CC,0}$.

Thus, \eqref{eq:uap} implies that
\begin{eqnarray*}
\RHDbmo_{\CC,0}&=& \sum_{-1\leq\reel\alpha<0}\sum_{p}\cD_{\CC,0}\cdot\cD_{\ov \CC,0}\cdot \wt u_{\alpha,p}
\end{eqnarray*}
from which we deduce that
\begin{eqnarray}
\RHDb_{\CC,0}&=& \sum_{-1\leq\reel\alpha<0}\sum_{p}\cD_{\CC,0}\cdot\cD_{\ov \CC,0}\cdot u_{\alpha,p}+\CC[\partial_t,\partial_{\ov t}]\cdot\delta\nonumber\\
&=&\sum_{-1\leq \reel\alpha<0}\sum_{p}\cD_{\CC,0}\cdot\cD_{\ov \CC,0}\cdot u_{\alpha,p}\label{eq:RHDb}
\end{eqnarray}
because the Dirac distribution $\delta$ can be written as
\begin{eqnarray*}
-2i\pi\delta&=&\partial_t\partial_{\ov t}\log\module{t}^2.
\end{eqnarray*}

\begin{Remarques}\label{rem:dvptasympt}
\begin{enumerate}
\item
Let $u\in\RHDb_{\CC,0}$ and put $\cM=\cD_{\CC,0}u\subset \RHDb_{\CC,0}$. Denote by $\wt u$ the image of $u$ in $\Dbmo_{\CC,0}$. Then the regular holonomic module $\cM[t^{-1}]$ is naturally embedded in $\RHDbmo_{\CC,0}$, the image of $u$ in $\cM[t^{-1}]$ is $\wt u$ \emph{via} this embedding and $\cM[t^{-1}]=\cD_{\CC,0}[t^{-1}]\wt u$, so $\cD_{\CC,0}\wt u$ is identified with the quotient of $\cM$ by its torsion supported at the origin.
\item
Using Borel's lemma, one can show that $\CRHDbmo_{\CC,0}$ is equal to the subspace of germs at of $C^\infty$ functions on $\CC^*$ having an infinitely termwise differentiable asymptotic expansion at $0$, in the sense of \cite{B-M87}, the exponents of which belong to a finite union of lattices in $\CC$.
\end{enumerate}
\end{Remarques}

\subsubsection*{Mellin transform}
We also have a characterization of $\CRHDb_{\CC,0}$ in terms of Mellin transform (\cite{B-M87}). Let $u\in\Db_{\CC,0}$ be a germ of distribution, and denote also by $u$ a representative of this germ in $\Gamma(D,\Db_\CC)$ where $D$ is a small disc centered at the origin. Let $\chi\in\cC_c^\infty(D)$ such that $\chi\equiv1$ near $0$ and having a sufficiently small support. Then, for any $k',k''\in\ZZ$, define
\begin{eqnarray*}
\cJ_{u}^{(k',k'')}(s)&=&\langle\chi u,t^{k'}\ov t^{k''}\module{t}^{2s}dt\wedge d\ov t\rangle
\end{eqnarray*}
which are holomorphic on $\reel(s)\gg0$. These functions depend on $\chi$ up to the addition of an entire function. So the classes of $\cJ_{u}^{(k',k'')}$ modulo $\cO(\CC)$ only depend on the germ $u$. Moreover, these functions can be recovered from the functions $\cJ_{u}^{(k,0)}$ and $\cJ_{u}^{(0,k)}$ for $k\in\NN$, because, if for instance $k'\geq k''$, we clearly have $\cJ_{u}^{(k',k'')}(s)=\cJ_{u}^{(k'-k'',0)}(s+k'')$. Moreover, $\cJ_{u}^{(k',k'')}$ only depends on the image of $u$ in $\Dbmo_{\CC,0}$.

\begin{proposition}[{\cite[Theorem 4]{B-M87}}]\label{prop:caractRH}
Let $u\in\Db_{\CC,0}$. Then $u\in \CRHDb_{\CC,0}$ if and only if there exists a finite set $\cR\subset \CC$ such that for all $k\in\NN$ the functions $\cJ_{u}^{(k,0)},\cJ_{u}^{(0,k)}$, which are holomorphic on $\reel(s)\gg0$, extend to meromorphic functions on $\CC$ with poles at most in $\cR+\ZZ$, and satisfy
\begin{multline*}
(\exists R>0),\ (\forall N>0),\ (\forall\ell>0),\ (\forall\ell'>0), \\ \module{s+k/2}^\ell\module{k}^{\ell'}\sup\lefpar \module{\cJ_{u}^{(k,0)}(s)},\module{\cJ_{u}^{(0,k)}(s)}\rigpar \leq C(u,N,\ell,\ell') R^{\reel(s+k/2)}
\end{multline*}
for $\reel(s+k/2+N)\geq -1$ and $\module{s+k/2}\gg0$.
\end{proposition}

\begin{proof}
Remark first that $u\in \RHDb_{\CC,0}$ if and only if its image in $\Dbmo_{\CC,0}$ belongs to $\RHDbmo_{\CC,0}$. Moreover, we may fix a representative for $u$ and consider $\chi u$ to define $\cJ_{u}^{(k,0)}(s)$ or $\cJ_{u}^{(0,k)}$. The condition in the proposition is easily seen to be independent of these choices. The result is then a direct consequence of \cite[Theorem 4]{B-M87}.
\end{proof}

\section{Hermitian duality and moderate nearby/vanishing cycles}

We will show in this section the compatibility between these functors. We will first recall briefly the construction of moderate and vanishing cycles for holonomic $\cD$-modules, in order to be able to give a detailed account of the compatibility.

\subsection{Notation}
We fix a total ordering on $\CC$, denoted by $\leq$, which is assumed to satisfy (a), (b), (c) below:
{\def\theenumi{\alph{enumi}}
\begin{enumerate}
\item
it induces the usual ordering on $\RR$,
\item
for $a\in\RR$, $\{z\in\CC\mid z<a\}=\{z\in\CC\mid \reel(z)<a\}$,
\item
for $a\in\RR$ and $z,z'\in\CC$, $z\leq z'\iff z+a\leq z'+a$.
\end{enumerate}}
In the following, we will choose the ordering on $\CC$ induced by lexicographically ordering the triples $(\reel(a),\module{\im(a)},\im(a))$. With such an ordering we have
\begin{eqnarray*}
\{\alpha\in\CC\mid -1\leq \alpha<0\}&=&\{\alpha\in\CC\mid -1\leq \reel\alpha<0\}.
\end{eqnarray*}
For a complex number $\gamma$, denote by $[\gamma]$ the largest integer less than or equal to $\gamma$, using the fixed total ordering on $\CC$.

\subsection{Review on the Malgrange-Kashiwara filtration}\label{subsec:vfil}
Let $Z$ be a complex analytic manifold of dimension $n$, put $X=Z\times \CC$, let $t$ denote the coordinate on $\CC$ or the projection $X\rightarrow \CC$ and consider the inclusion $Z=Z\times \{0\}\hookrightarrow X$.

For a holonomic $\cD_X$-module $\cM$, let $V_\bbullet(\cM)$ be the Malgrange-Kashiwara filtration on $\cM$ relative to $Z\times \{0\}$ (see \eg. \cite{M-S86}): this is a filtration indexed by the union of a (locally on $Z$) finite number of lattices $\sigma+\ZZ\subset\CC$ ($\sigma\in\ccS$ and we may choose the finite set $\ccS\subset\CC$ contained in $\reel(\sigma)\in[0,1[$), using the ordering specified above. For any $\alpha\in\CC$, the graded module $\gr_{\alpha}^{V}\cM\defin V_\alpha\cM/V_{<\alpha}\cM$ is $\cD_Z$-holonomic (and moreover regular when $\cM$ is so) and comes equipped with a nilpotent endomorphism $N$, induced by the action of $-(\partial_tt+\alpha)$.

We have isomorphisms
\begin{equation}\label{eq:isot}
t:V_\alpha\cM\isom V_{\alpha-1}\cM\quad (\alpha<0)
\end{equation}
and
\begin{equation}\label{eq:isodt}
\partial_t:\gr_{\alpha}^{V}\cM\isom\gr_{\alpha+1}^{V}\cM\quad (\alpha>-1).
\end{equation}
The complex $i^+\cM$ is quasi-isomorphic to the complex
\[
\gr_{0}^{V}\cM\MRE{t}\gr_{-1}^{V}\cM
\]
(where the right term has degree $0$) and if $\cM=j_\dag j^+\cM$ it is also isomorphic to the complex
\[
\gr_{-1}^{V}\cM\MRE{t\partial_t}\gr_{-1}^{V}\cM.
\]
Similarly, the complex $i^\dag\cM$ is quasi-isomorphic to the complex
\[
\gr_{-1}^{V}\cM\MRE{\partial_t}\gr_{0}^{V}\cM
\]
(where the left term has degree $0$) and if $\cM=j_+ j^+\cM$ it is also isomorphic to the complex
\[
\gr_{-1}^{V}\cM\MRE{t\partial_t}\gr_{-1}^{V}\cM.
\]
In particular, if $\cM=j_+ j^+\cM$, we will identify
\begin{align}\label{eq:idag}
\cH^0(i^\dag\cM)&\qqbox{with} \ker[t\partial_t:\gr_{-1}^{V}\cM\rightarrow \gr_{-1}^{V}\cM]
\end{align}
and, if $\cM=j_\dag j^+\cM$,
\begin{align}\label{eq:iplus}
\cH^0(i^+\cM)&\qqbox{with} \coker[t\partial_t:\gr_{-1}^{V}\cM\rightarrow \gr_{-1}^{V}\cM].
\end{align}

Analogous results hold for holonomic $\cD_{\ov X}$-modules. We still denote by $V_\bbullet$ the Mal\-gran\-ge-Kashiwa\-ra filtration and by $N$ the nilpotent endomorphism induced by $-(\partial_{\ov t}\ov t+\alpha)$ on $\gr_{\alpha}^{V}$.

\subsection{Review on moderate nearby and vanishing cycles (see \eg.\ \cite{M-S86,MSaito86})} \label{subsec:nearby}
Let $\cM$ be a holonomic $\cD_X$-module (specializable would be enough, see \eg.\ \cite{M-S86}). Let $\alpha$ be such that $-1\leq \alpha<0$ and put $\lambda=\exp(2i\pi\alpha)$. For $p\in\NN$, put $\cM_{\alpha,p}=(\cM[t^{-1}])^{p+1}=\oplus_{k=0}^{p}\cM[t^{-1}]\otimes e_{\alpha,k}$. The $\cD_{X/\CC}$-structure on $\cM_{\alpha,p}$ is the direct sum of the $\cD_{X/\CC}$-structures on each term $\cM[t^{-1}]$ and the $\cD_X$-structure is given by the relation
\[
t\partial_t(m\otimes e_{\alpha,k})=[(\partial_tt+\alpha)m]\otimes e_{\alpha,k}+m\otimes e_{\alpha,k-1},
\]
with the convention that $e_{\alpha,k}=0$ for $k<0$. Remark that $\cM[t^{-1}]$ is a direct summand of $\cM_{-1,p}$ for any $p\geq 0$ (we may consider that $e_{\alpha,k}$ plays the role of the multivalued function $t^{\alpha+1}(\log t)^k/k!$).

We have natural morphisms of $\cD_X$-modules:
\[
\begin{array}{rcl}
\cM_{\alpha,p}&\HMRE{a_{p,p+1}}&\cM_{\alpha,p+1}\\
\dpl\sum_{k=0}^{p}m_{\alpha,k}\otimes e_{\alpha,k}&\DMRE{}&\dpl\sum_{k=0}^{p}m_{\alpha,k}\otimes e_{\alpha,k}
\end{array}
\]
and
\[
\begin{array}{rcl}
\cM_{\alpha,p+1}&\MRE{b_{p+1,p}}&\cM_{\alpha,p}\\
\dpl\sum_{k=0}^{p+1}m_{\alpha,k}\otimes e_{\alpha,k}&\DMRE{}&\dpl\sum_{k=0}^{p}m_{\alpha,k+1}\otimes e_{\alpha,k}.
\end{array}
\]
We will denote by $N$ (without index $p$) any of the endomorphisms
\[
N=a_{p-1,p}\circ b_{p,p-1}:\cM_{\alpha,p}\longrightarrow \cM_{\alpha,p},
\]
sending $m\otimes e_{\alpha,k}$ to $m\otimes e_{\alpha,k-1}$.
The inductive (\resp.\ projective) system $\cH^0(i^\dag\cM_{\alpha,p})$ (\resp.\ $\cH^0(i^+j_\dag j^+\cM_{\alpha,p})$) where the maps are induced by $a_{p,p+1}$ (\resp.\ $b_{p+1,p}$) is stationary locally on $X$, and both systems have a common limit isomorphic to $\gr_{\alpha}^{V}\cM$: we may identify $\gr_{-1}^{V}\cM_{\alpha,p}$ with $\oplus_{k=0}^{p}\gr_{\alpha}^{V}\cM\otimes e_{\alpha,k}$; the natural mappings
\begin{align}
\gr_{\alpha}^{V}\cM&\MRE{} \gr_{-1}^{V}\cM_{\alpha,p}\notag\\
m_0&\DMRE{}\ooplus_{k=0}^{p}[-(\partial_tt+\alpha)]^km_0\otimes e_{\alpha,k}
\label{eq:a0p}\\
\tag*{and}\\
\gr_{-1}^{V}\cM_{\alpha,p}&\MRE{} \gr_{\alpha}^{V}\cM\notag\\
\ooplus_{k=0}^{p}m_k\otimes e_{\alpha,k}&\DMRE{}\ooplus_{k=0}^{p}[-(\partial_tt+\alpha)]^km_{p-k}\label{eq:bp0}
\end{align}
induce, for $p$ large enough, an isomorphism from $\gr_{\alpha}^{V}\cM$ to $\ker t\partial_t\simeq\cH^0(i^\dag\cM_{\alpha,p})$ and from $\coker t\partial_t\simeq \cH^0(i^+j_\dag j^+\cM_{\alpha,p})$ to $\gr_{\alpha}^{V}\cM$.

We denote this limit by $\psim_{t,\lambda}\cM$ and call it the moderate nearby cycle module associated with $\cM$, with eigenvalue $\lambda$. We also denote by $N$ the endomorphism induced by the previous $N$. It corresponds naturally to $-(\partial_tt+\alpha)$ \emph{via} both isomorphisms with $\gr_{\alpha}^{V}\cM$. Notice also that the inductive system of $\cH^1$ (\resp.\ the projective system of $\cH^{-1}$) has limit $0$.

\medskip
The construction of the moderate vanishing cycle module $\phim_{t,1}(\cM)$ is achieved by considering the inductive sytem of complexes $\cM\rightarrow \cM_{-1,p}$ (where the right term has degree $0$ and the map is the composition of $\loc:\cM\rightarrow \cM[t^{-1}]$ with $a_{0,p}: \cM[t^{-1}]\rightarrow \cM_{-1,p}$) instead of the single module $\cM_{\alpha,p}$. The only possible non vanishing limit is also obtained for $\cH^0i^\dag$. It can also be achieved by considering the projective system of complexes $j_\dag j^+\cM_{-1,p}\rightarrow \cM$ (where the left term has degree $0$ and the map is the composition of $j_\dag j^+ b_{p,0}$ and $\coloc: j_\dag j^+\cM\rightarrow \cM$) and the projective limit of $\cH^0i^+$. Let us give some precise description. The complex $i^\dag(\cM\rightarrow \cM_{-1,p})$ is the single complex associated to the double complex
\[
\begin{array}{c}
\xymatrix{
j_\dag j^+\cM\ar[r]\ar[d]_-{\coloc}&j_\dag j^+\cM_{-1,p}\ar[d]^-{\coloc}\\
\cM\ar[r]&\cM_{-1,p}
}
\end{array}
\quad\simeq\quad
\begin{array}{c}
\xymatrix{
\gr_{-1}^{V}\cM\ar[r]\ar[d]_-{\partial_t}& \gr_{-1}^{V}(\cM_{-1,p})\ar[d]^-{t\partial_t}\\
\gr_0^V\cM\ar[r]_-t&\gr_{-1}^{V}(\cM_{-1,p})
}
\end{array}
\]
which is isomorphic to the complex
\[
\gr_{-1}^{V}\cM\longrightarrow \gr_0^V\cM\oplus \gr_{-1}^{V}(\cM_{-1,p})\longrightarrow \gr_{-1}^{V}(\cM_{-1,p})
\]
where the middle term has degree $0$. The kernel of the second morphism can be identified with $\gr_{-1}^{V}\cM\oplus \gr_{0}^{V}\cM$ \emph{via}
\[
m_0\oplus n_0\longmapsto n_0\oplus (m_0\otimes e_{-1,0})\oplus \lefcro\ooplus_{k=1}^{p}(-t\partial_t)^{k-1}(-t\partial_tm_0+tn_0)\otimes e_{-1,k}\rigcro
\]
and the $\cH^0$ of this complex is identified to $\gr_0^V\cM$ \emph{via}
\[
\gr_0^V\cM\MRE{0\oplus \id}\gr_{-1}^{V}\cM\oplus \gr_{0}^{V}\cM.
\]
The action of $0\oplus N$ on $\gr_0^V\cM\oplus \gr_{-1}^{V}(\cM_{-1,p})$ induces, \emph{via} these isomorphisms, the action of $-\partial_tt$ on $\gr_0^V\cM$.

Similarly, the complex $i^+(j_\dag j^+\cM_{-1,p}\rightarrow \cM)$ is isomorphic to the single complex associated with
\[
\xymatrix{
\gr_{-1}^{V}(\cM_{-1,p})\ar[d]_-{t\partial_t}\ar[r]^-{\partial_t}& \gr_0^V\cM\ar[d]^-t\\
\gr_{-1}^{V}(\cM_{-1,p})\ar[r]&\gr_{-1}^{V}\cM
}
\]
where the middle term has degree $0$. Its $\cH^0$ is naturally
isomorphic to $\gr_0^V\cM$ and the action of $N$ on $\cM_{-1,p}$
induces that of $-\partial_tt$ on $\gr_0^V\cM$.

\medskip
The morphisms $\can$ and $\Var$ are defined as
\[
\xymatrix @=2cm{
\gr_{-1}^{V}\cM\ar@/^/[r]^-{\can=-\partial_t}& \gr_{0}^{V}\cM\ar@/^/[l]^-{\Var=t}
}
\]
and can be obtained, \emph{via} the previous isomorphisms, as coming
from the morphisms of complexes
$$
\begin{array}{c}
\xymatrix{
0\ar[r]\ar[d]&\cM_{-1,p}\ar[d]^-{\id}\\
\cM\ar[r]&\cM_{-1,p}
}
\end{array}
\qqbox{or}
\begin{array}{c}
\xymatrix{
j_\dag j^+\cM_{-1,p}\ar[d]_-{N}\ar[r]\ar[d]&0\ar[d]\\
j_\dag j^+\cM_{-1,p}\ar[r]&\cM
}
\end{array}
\leqno{(\can)}
$$
$$
\begin{array}{c}
\xymatrix{
\cM\ar[r]\ar[d]&\cM_{-1,p}\ar[d]^-{N}\\
0\ar[r]&\cM_{-1,p}
}
\end{array}
\qqbox{or}
\begin{array}{c}
\xymatrix{
j_\dag j^+\cM_{-1,p}\ar[d]_-{\id}\ar[r]\ar[d]&\cM\ar[d]\\
j_\dag j^+\cM_{-1,p}\ar[r]&0
}
\end{array}
\leqno{(\Var)}
$$

\subsection{Compatibility with Hermitian duality}
We now assume that $\cM$ is regular holonomic. For any $\alpha$ such
that $-1\leq \alpha<0$, consider the function $u_{-\alpha-2,p}=\module{t}^{2(\alpha+1)}(\log\module{t}^2)^k/k!$
analogous to that of \T\ref{subsec:RH1} as a function on $X$. It has
moderate growth along $Z$ as well as all its derivatives. Hence, for
any moderate distribution $\wt u$ along $Z$, the product $u_{-\alpha-2,p}\wt u$ is well-defined as a moderate distribution along $Z$.

\begin{lemme}\label{lem:herm}
The pairing
\begin{eqnarray*}
C_X(\cM)_{\alpha,p}\ootimes_\CC\cM_{\alpha,p}&\longrightarrow &\Dbmz_X\\
\lefpar \sum_{k=0}^{p}\mu_{\alpha,k}\otimes \ov e_{\alpha,k}\rigpar \otimes \lefpar \sum_{\ell=0}^{p}m_{\alpha,\ell}\otimes e_{\alpha,\ell}\rigpar &\longmapsto&\sum_{k,\ell}\mu_{\alpha,k}(m_{\alpha,\ell})u_{-\alpha,k+\ell-p}
\end{eqnarray*}
is nondegenerate and induces an isomorphism compatible with $N$ and $\Cmz_X(N)$
\[
 \eta_{\alpha,p}:C_X(\cM)_{\alpha,p}\isom\Cmz_X(\cM_{\alpha,p})
\]
such that all diagrams
\[
\begin{array}{c}
\xymatrix{
C_X(\cM)_{\alpha,p}\ar[d]_-{a_{p,p+1}}\ar[r]^-{\sim}&\Cmz_X(\cM_{\alpha,p})\ar[d]^-{\Cmz_X(b_{p+1,p})}\\
C_X(\cM)_{\alpha,p+1}\ar[r]^-{\sim}&\Cmz_X(\cM_{\alpha,p+1})
}
\end{array}
\]
and
\[
\begin{array}{c}
\xymatrix{
C_X(\cM)_{\alpha,p}\ar[r]^-{\sim}&\Cmz_X(\cM_{\alpha,p})\\
C_X(\cM)_{\alpha,p+1}\ar[u]^-{b_{p+1,p}}\ar[r]^-{\sim}&\Cmz_X(\cM_{\alpha,p+1})\ar[u]_-{\Cmz_X(a_{p,p+1})}
}
\end{array}
\]
commute.
\end{lemme}

\begin{proof}
First, it is easy to see that the morphism $\eta_{\alpha,p}$ induced by the pairing induces commutative diagrams as in the lemma. The compatibility of $\eta_{\alpha,p}$ with $N$ and $\Cmz_X(N)$ is thus clear. The nondegeneracy of the pairing is then proved by induction on $p$, the case $p=0$ being easy.
\end{proof}

\begin{theoreme}\label{th:compatible}
There exist natural isomorphisms of functors from $\Mod_{hr}(\cD_X)$ to $\Mod_{hr}(\cD_{\ov Z})$
\[
c_{X,\lambda}^{\psi}:\psim_{t,\lambda}\circ C_X\longrightarrow C_Z\circ \psim_{t,\lambda}, \quad(\lambda\in\CC^*)\qqbox{and}c_{X,1}^{\phi}:\phim_{t,1}\circ C_X\longrightarrow C_Z\circ \phim_{t,1}
\]
which satisfy the following properties, putting $c_X=c_{X,\lambda}^{\psi}$ or $c_{X,1}^{\phi}$:
\begin{itemize}
\item
$c_X=C_Z\circ c_{\ov X}\circ C_X$;
\item
$c_X\circ N=C_Z(N)\circ c_X$;
\item
$c_{X,1}^{\phi}\circ \can=C_Z(\Var)\circ c_{X,1}^{\psi}$ and $c_{X,1}^{\psi}\circ \Var=C_Z(\can)\circ c_{X,1}^{\phi}$.
\end{itemize}
\end{theoreme}

\begin{proof}
According to the previous lemma and to Corollary \ref{cor:imdir}, the inductive system $$(\cH^0\ov\imath^\dag C_X(\cM)_{\alpha,p},\cH^0\ov\imath^\dag a_{p,p+1})$$ is isomorphic, \emph{via} $\cH^0\ov\imath^\dag \eta_{\alpha,p}$, to $C_Z$ of the projective system $(\cH^0(i^+j_\dag j^+\cM_{\alpha,p}), b_{p+1,p})$. The first part of the theorem then follows from the construction of $\psim_{t,\lambda}$ recalled in \T\ref{subsec:nearby}. The proof for $\phim_{t,1}$ and the other properties also follow from the same arguments.
\end{proof}

\subsection{Nearby/vanishing cycles for a sesquilinear form} \label{subsec:sesqui}
Let $\cM',\cM''$ be two regular holonomic $\cD_X$-modules and let $S:\cM'\otimes_\CC\ov{\cM''}\rightarrow \Db_X$ be a sesquilinear pairing.

We will define, for $-1\leq \alpha< 0$, sesquilinear forms
\[
\psi_\lambda S: \gr_{\alpha}^{V}\cM'\ootimes_{\CC}\gr_{\alpha}^{V}\ov{\cM''}\longrightarrow \Db_Z
\]
and similarly (for $\alpha=0$) $\phi_1 S$, which satisfy (with obvious notation)
\begin{equation}\label{eq:propS}
\begin{split}
\psi_\lambda S(N\cbbullet,\cbbullet)&=\psi_\lambda S(\cbbullet,N\cbbullet)\\
\phi_1 S(N\cbbullet,\cbbullet)&=\phi_1 S(\cbbullet,N\cbbullet)\\
\psi_{1}S(\Var\cbbullet,\cbbullet)&=\phi_1 S(\cbbullet,\can\cbbullet)\\
\psi_{1}S(\cbbullet,\Var\cbbullet)&=\phi_1 S(\can\cbbullet,\cbbullet).
\end{split}
\end{equation}

Denote for a while by $L_S$ the $\cD_{\ov X}$-linear morphism $\ov{\cM''}\rightarrow C_X\cM'$ induced by $S$. Consider $\psi_\lambda L_S:\gr_\alpha^V\ov{\cM''}\rightarrow \gr_\alpha^VC_X\cM'$ (and $\phi_1L_S$ defined similarly). Its composition with $c_{X,\lambda}^{\psi}$ (or $c_{X,1}^{\phi}$) is the linear morphism associated with a sesquilinear form $\psi_\lambda S$ or $\phi_1 S$. The properties \eqref{eq:propS} follow then from the properties of $c_X$ given by Theorem \ref{th:compatible}.

\begin{remarque}\label{rem:prim}
Denote by $M_\bbullet\gr_\alpha^V(\cM)$ the monodromy filtration associated to the nilpotent endomorphism $N$, \ie.\ the increasing filtration such that $NM_k\subset M_{k-2}$ and for all $\ell\geq 0$,
\begin{eqnarray*}
\gr_{\ell}^{M}\gr_\alpha^V\cM&\MRE{N^\ell}&\gr_{-\ell}^{M}\gr_\alpha^V\cM
\end{eqnarray*}
is an isomorphism. Let $P\gr_{\ell}^{M}\gr_\alpha^V\cM$ denote the primitive part
\[
\ker\lefcro N^{\ell+1}:\gr_{\ell}^{M}\gr_\alpha^V\cM\longrightarrow \gr_{-\ell-2}^{M}\gr_\alpha^V\cM\rigcro.
\]
The pairing $\psi_\lambda S$, being compatible with $N$, induces for any $\ell$ a pairing
\begin{eqnarray*}
\gr_{\ell}^{M}\gr_\alpha^V\cM'\ootimes_\CC\gr_{-\ell}^{M}\gr_\alpha^V\ov{\cM''}&\MRE{\psi_{\lambda,\ell}S}&\Db_Z
\end{eqnarray*}
and is nondegenerate iff $\psi_{\lambda,\ell}S$ is nondegenerate for any $\ell$. This is so iff the pairing induced on the primitive parts
\begin{eqnarray}\label{eq:conjprim}
P\gr_{\ell}^{M}\gr_\alpha^V\cM'\ootimes_\CC P\gr_{\ell}^{M}\gr_\alpha^V\ov{\cM''}&\MRE{\psi_{\lambda,\ell}S\circ(\id\otimes N^\ell)}&\Db_Z
\end{eqnarray}
is nondegenerate, according to the Lefschetz decomposition. Similar results hold for $\phi_1S$.

For $\ell\geq 0$ we will set
\[
P\psi_{\lambda,\ell}S\defin\psi_{\lambda,\ell}S\circ(\id\otimes N^\ell)\qqbox{and} P\phi_{1,\ell}S\defin\phi_{1,\ell}S\circ(\id\otimes N^\ell)
\]
\end{remarque}

We deduce from Theorem \ref{th:compatible}:

\begin{corollaire}\label{cor:sesqui}
The sesquilinear form $S$ is nondegenerate in a neighbourhood of $Z$ if and only if all sesquilinear forms $P\psi_{\lambda,\ell}S$ ($\lambda\in\CC^*$, $\ell\geq 0$) and $P\phi_{1,\ell}S$ ($\ell\geq 0$) are nondegenerate.
\end{corollaire}

\begin{proof}
According to Remark \ref{rem:prim}, it is enough to show that $S$ is nondegenerate iff all $\psi_\lambda S$ and $\phi_1S$ are so. 
Now, $L_S$ is an isomorphism in a neighbourhood of $Z$ if and only if all $\psi_\lambda L_S$ and $\phi_1L_S$ are isomorphisms: this follows from the fact that a regular holonomic module $\cM$ is equal to zero near $Z$ if and only if all its moderate nearby or vanishing cycles vanish on $Z$. The result is then a consequence of the definition of $\psi_\lambda S$ and $\phi_1S$ and of Theorem \ref{th:compatible}.
\end{proof}

\section{Hermitian duality and asymptotic expansions}
We will give in this section a more explicit description of the compatibility morphisms given in Theorem \ref{th:compatible}, using asymptotic expansions (in the sense of Remark \ref{rem:dvptasympt}(2)). The main goal will be to give a more precise version of Theorem \ref{th:bm}, taking into account the order with respect to the Malgrange-Kashiwara filtration.

We begin with some easy results in dimension $1$.

\subsection{Dimension $1$}
\subsubsection*{Regular holonomic distributions and Malgrange-Kashiwara filtration}
If $u$ is the germ at $0\in\CC$ of a regular holonomic distribution defined on some open disc $D$ centered at the origin, we denote by $\alpha'(u)$ the order of $u$ with respect to the Malgrange-Kashiwara filtration of the regular holonomic module $\cD_Du\subset\Db_D$ and by $\alpha''(u)$ its $V_\bbullet$-order in $\cD_{\ov D}u$. Notice that, according to the strictness property of any morphism between holonomic modules with respect to the Malgrange-Kashiwara filtration, if $v\in\cD_D\cdot u$, then $\alpha'(v)$ is equal to the $V$-order of $v$ when viewed as an element of $\cD_D\cdot u$.

We obtain in this way increasing filtrations
\begin{eqnarray*}
V'_{\alpha'}(\RHDb_{\CC,0})&=&\{u\in\RHDb_{\CC,0}\mid \alpha'(u)\leq \alpha'\}\\
V''_{\alpha''}(\RHDb_{\CC,0})&=&\{u\in\RHDb_{\CC,0}\mid \alpha''(u)\leq \alpha''\}
\end{eqnarray*}
(where $\leq$ is the fixed total ordering on $\CC$) and thus a doubly indexed filtration
\begin{eqnarray*}
V_{\alpha',\alpha''}(\RHDb_{\CC,0})&=&V'_{\alpha'}(\RHDb_{\CC,0})\cap V''_{\alpha''}(\RHDb_{\CC,0}).
\end{eqnarray*}
We then put
\begin{eqnarray*}
\gr_{\alpha',\alpha''}^{V}\RHDb_{\CC,0}&\defin&V_{\alpha',\alpha''}(\RHDb_{\CC,0})/(V_{<\alpha',\alpha''}+V_{\alpha',<\alpha''})(\RHDb_{\CC,0}).
\end{eqnarray*}

For $\lambda\in\CC^*$, choose $\alpha\in\CC$ with $-1\leq \alpha<0$ such that $\lambda=\exp(2i\pi\alpha)$; put as in \T\ref{subsec:RH1}, $u_{\alpha,p}=\module{t}^{-2(\alpha+1)}\dfrac{(\log\module{t}^2)^p}{p!}$ and
\begin{equation*}
\RHDb_{\CC,0}(\lambda)=
\sum_{p\geq 0}\cD_{\CC,0}\cD_{\ov\CC,0}\cdot u_{\alpha,p}
\end{equation*}
We then have
\begin{eqnarray*}
\RHDb_{\CC,0}&=&\ooplus_{\lambda\in\CC^*}\RHDb_{\CC,0}(\lambda).
\end{eqnarray*}

\begin{proposition}\label{prop:VRH}
The filtration $V_{\bbullet,\bbullet}(\RHDb_{\CC,0})$ satisfies the following properties.
\begin{enumerate}
\item\label{VRH-1}
$tV_{\alpha',\alpha''}\subset V_{\alpha'-1,\alpha''}$, \resp.\ $\ov t V_{\alpha',\alpha''}\subset V_{\alpha',\alpha''-1}$, with equality if $\alpha'<0$, \resp.\ $\alpha''<0$.
\item\label{VRH-2}
$\partial_tV_{\alpha',\alpha''}\subset V_{\alpha'+1,\alpha''}$, \resp.\ $\partial_{\ov t}V_{\alpha',\alpha''}\subset V_{\alpha',\alpha''+1}$.
\item\label{VRH-3}
Let $u\in\RHDb_{\CC,0}$. Then $u\in V_{\alpha',\alpha''}$ iff there exist $k',k''\in\NN$ with
\[
(\partial_tt+\alpha')^{k'}u\in V_{<\alpha',\alpha''} \quad\text{and}\quad (\partial_{\ov t}\ov t+\alpha'')^{k''}u\in V_{\alpha',<\alpha''}.
\]
\item\label{VRH-4}
We have
\begin{eqnarray*}
V_{\alpha',\alpha''}(\RHDb_{\CC,0})&=&\sum_{-1\leq \alpha<0}\sum_{\substack{k',k''\in\ZZ\\ \alpha+k'\leq \alpha'\\ \alpha+k''\leq \alpha''}}\sum_{p\geq 0} V_{k'}(\cD_{\CC,0})\cdot V_{k''}(\cD_{\ov\CC,0})\cdot u_{\alpha,p}
\end{eqnarray*}
and $\gr_{\alpha',\alpha''}^{V}(\RHDb_{\CC,0})=0$ if $\alpha'-\alpha''\not\in\ZZ$.
\item\label{VRH-5}
For any $\alpha\in\CC$, $(\partial_tt-\partial_{\ov t}\ov t)V_{\alpha,\alpha}\subset V_{<\alpha,\alpha}+V_{\alpha,<\alpha}$.
\end{enumerate}
\end{proposition}

\begin{proof}
The assertion (\ref{VRH-2}) and the first part of (\ref{VRH-1}) follow immediately from the properties of the Malgrange-Kashiwara filtration on holonomic modules.

Let us prove the second part of (\ref{VRH-1}). Let $u\in V_{\alpha'-1,\alpha''}(\RHDb_{\CC,0})$ with $\alpha'<0$. According to \eqref{eq:isot}, there exists then $v\in\RHDb_{\CC,0}$ with $\alpha'(v)\leq \alpha'$ such that $u=tv$. Let $b(\partial_{\ov t}\ov t)$ be the minimal polynomial satisfying $b(\partial_{\ov t}\ov t)u=\ov t P(\ov t,\partial_{\ov t}\ov t)u$ with $P\in V_0\cD_{\ov D,0}$. Then $w\defin [b(\partial_{\ov t}\ov t)-\ov t P(\ov t,\partial_{\ov t}\ov t)]v$ is supported at the origin and satisfies $\alpha'(w)\leq \alpha'(v)\leq \alpha'<0$. Therefore, by \eqref{eq:isot}, we have $w=0$ and $v$ satisfies $\alpha''(v)\leq \alpha''(u)$.

\smallskip
For (\ref{VRH-3}), remark that there exist $\beta',\beta''$ such that $u\in V_{\beta',\beta''}(\RHDb_{\CC,0})$. We may assume that $\beta'\geq \alpha'$ and $\beta''\geq \alpha''$. There exists polynomials $B'(-s)$ (\resp.\ $B''(-s)$) with roots in $]\alpha',\beta']$ (\resp.\ in $]\alpha'',\beta'']$), such that $B'(\partial_tt)B''(\partial_{\ov t}\ov t)u$ belongs to $V_{\alpha',\alpha''}$. Applying B\'ezout and the condition in (\ref{VRH-3}) we conclude that $u$ belongs to $V_{\alpha',\alpha''}$.

\smallskip
Let us now prove (\ref{VRH-4}) and (\ref{VRH-5}). We will first need the following lemma.
\begin{lemme}\label{lem:VaaRH}
{\def\theenumi{\alph{enumi}}
\begin{enumerate}
\item
We have $\gr^V_{\beta',\beta''}\RHDb_{\CC,0}(\lambda)=0$ if $\beta'\not\in\alpha+\ZZ$ or $\beta''\not\in\alpha+\ZZ$.
\item
For all $k',k''\in\ZZ$ we have
\begin{equation*}
V_{k'+\alpha,k''+\alpha}(\RHDb_{\CC,0}(\lambda))=
\sum_{p\geq 0}V_{k'}(\cD_{\CC,0})V_{k''}(\cD_{\ov\CC,0})\cdot u_{\alpha,p}.
\end{equation*}
\item
For $-1\leq \alpha<0$, the classes of $u_{\alpha,p}$ ($p\geq 0$) form a basis of the $\CC$-vector space $\gr_{\alpha,\alpha}^{V}\RHDb_{\CC,0}$.
\item
The classes of $\partial_t\partial_{\ov t}u_{-1,p}$ ($p\geq 1$) form a basis of $\gr_{0,0}^{V}\RHDb_{\CC,0}$.
\end{enumerate}
}
\end{lemme}

\begin{proof}
According to \eqref{eq:uap}, the distribution $u_{\alpha,p}$ ($-1\leq \alpha<0$ and $p\in\NN$) satisfies
\begin{equation}\label{eq:uapb}
(\partial_tt+\alpha)^{p+1}u_{\alpha,p}=(\partial_{\ov t}\ov t+\alpha)^{p+1}u_{\alpha,p}=0.
\end{equation}
It is then in $V_{\alpha,\alpha}$.

It follows that, for any $\ov P\in\cD_{\ov \CC,0}$, the correspondence $1\mapsto \ov P u_{\alpha,p}$ induces a surjective $\cD_{\CC,0}$-linear morphism
\begin{eqnarray*}
\cD_{\CC,0}/\cD_{\CC,0}(\partial_tt+\alpha)^{p+1}&\longrightarrow &\cD_{\CC,0}\cdot \ov P u_{\alpha,p}.
\end{eqnarray*}
This implies that, for any $k\in\ZZ$, we have
\begin{eqnarray}\label{eq:Vka}
V'_{k+\alpha}\lefpar \cD_{\CC,0}\cdot \ov P u_{\alpha,p}\rigpar &=&V_k(\cD_{\CC,0})\cdot \ov P u_{\alpha,p}
\end{eqnarray}
because a similar property is easily seen to be true for $\cD_{\CC,0}/\cD_{\CC,0}(\partial_tt+\alpha)^{p+1}$ and any morphism of holonomic $\cD$-modules is strict with respect to the Malgrange-Kashiwara filtration.

By the same argument we also get that, for $-1\leq \alpha<0$,
\begin{eqnarray*}
V_{k'+\alpha,k''+\alpha}\lefpar \cD_{\CC,0}\cD_{\ov\CC,0}u_{\alpha,p}\rigpar &=&V_{k'}(\cD_{\CC,0})V_{k''}(\cD_{\ov\CC,0})u_{\alpha,p}.
\end{eqnarray*}
As $\RHDb_{\CC,0}(\lambda)=\varinjlim_p\cD_{\CC,0}\cD_{\ov\CC,0}u_{\alpha,p}$, the statements (a) and (b) are clear.

\medskip
Part (b) shows that the elements given in part (c) or (d) generate the corresponding bigraded object. If we have, for $-1\leq \alpha<0$, a linear relation between the classes of $u_{\alpha,p}$ ($p\geq 0$) in $\gr_{\alpha,\alpha}^{V}\RHDb_{\CC,0}$ then, by applying a suitable power of $\partial_tt+\alpha$ and using relation \eqref{eq:uap}, we would have a relation
\begin{eqnarray*}
\module{t}^{-2(\alpha+1)}&\in&t\sum_p\cO_{\CC,0}\cO_{\ov\CC,0}\cdot u_{\alpha,p}+\ov t\sum_p\cO_{\CC,0}\cO_{\ov\CC,0}\cdot u_{\alpha,p},
\end{eqnarray*}
which is clearly impossible by considering the valuation at $0$. Similarly, a linear relation between the classes $\partial_t\partial_{\ov t}u_{-1,p}$ ($p\geq 1$) would imply that $\delta\in (V_{-1,0}+V_{0,-1})\RHDb_{\CC,0}$. Notice now that
\[
t:V_{-1,0}(\RHDb_{\CC,0})\longrightarrow V_{-2,0}(\RHDb_{\CC,0})
\] 
is bijective: part (b) shows that it is onto; it is injective because $t:V'_{-1}\rightarrow V'_{-2}$ is so, as follows from \eqref{eq:isot}.

So, if $\delta=u^{(-1,0)}+u^{(0,-1)}$, we have $tu^{(-1,0)}\in V_{-2,0}\cap V_{-1,-1}=V_{-2,-1}$, hence $u^{(-1,0)}\in V_{-1,-1}$ and similarly $u^{(0,-1)}\in V_{-1,-1}$, so $\delta\in V_{-1,-1}$, which is impossible because $t$ acting on $V_{-1,-1}$ is injective.
\end{proof}

The statement (\ref{VRH-4}) of the proposition follows from (b) in the lemma.

\eqref{eq:uap} clearly implies that, for all $k,\ell\in\NN$, we have $(\partial_tt-\partial_{\ov t}\ov t)u=0$ if $u=t^k\ov t^k\partial_{t}^{\ell}\partial_{\ov t}^{\ell}u_{\alpha,p}$, for $-1\leq \alpha<0$. Then (\ref{VRH-5}) follows immediately.
\end{proof}

\subsubsection*{The Malgrange-Kashiwara filtration for $\CRHDb_{\CC,0}$}
In order to apply similar considerations to asymptotic expansion, we will introduce the Malgrange-Kashiwara filtration on $\CRHDb_{\CC,0}$. Put
\begin{eqnarray*}
V_{\alpha',\alpha''}(\CRHDb_{\CC,0})&\defin&\cC^\infty \cdot V_{\alpha',\alpha''}(\RHDb_{\CC,0}).
\end{eqnarray*}
We clearly have
\begin{align}\label{eq:Vaap}
V_{\alpha',\alpha''}(\CRHDb_{\CC,0})&=V_{\alpha',\alpha''}(\RHDb_{\CC,0})+(V_{<\alpha',\alpha''}+V_{\alpha',<\alpha''})(\CRHDb_{\CC,0}),
\end{align}
hence a surjective morphism $\gr_{\alpha',\alpha''}^{V}\RHDb_{\CC,0}\rightarrow \gr_{\alpha',\alpha''}^{V}\CRHDb_{\CC,0}$.

\begin{proposition}\label{prop:VCRH}
The results of Proposition {\rm \ref{prop:VRH}} apply as well to $V_{\bbullet,\bbullet}(\CRHDb_{\CC,0})$ and moreover
\begin{eqnarray*}
\gr_{\alpha',\alpha''}^{V}\RHDb_{\CC,0}=\gr_{\alpha',\alpha''}^{V}\CRHDb_{\CC,0}.
\end{eqnarray*}
\end{proposition}

\begin{remarque}\label{rem:from5}
It follows from Propositions \ref{prop:VRH}(\ref{VRH-5}) and \ref{prop:VCRH} that the nilpotent endomorphisms induced by $\partial_tt+\alpha$ or $\partial_{\ov t}\ov t+\alpha$ on $\gr_{\alpha,\alpha}^{V}\CRHDb_{\CC,0}$ coincide, for $-1\leq \alpha\leq 0$.
\end{remarque}

\begin{proof}
(\ref{VRH-1}), (\ref{VRH-2}), (\ref{VRH-3}) and (\ref{VRH-5}) in \ref{prop:VRH} immediately extend to $\CRHDb_{\CC,0}$. Moreover, (\ref{VRH-4}) clearly gives
\begin{eqnarray}\label{eq:VCRH-4}
V_{\alpha',\alpha''}(\CRHDb_{\CC,0})&=&\sum_{-1\leq \alpha<0}\sum_{\substack{k',k''\in\ZZ\\ \alpha+k'\leq \alpha'\\ \alpha+k''\leq \alpha''}}\sum_{p\geq 0} \cC^\infty\cdot V_{k'}(\cD_{\CC,0})\cdot V_{k''}(\cD_{\ov\CC,0})\cdot u_{\alpha,p}.
\end{eqnarray}
An argument similar to that of (c) and (d) in Lemma \ref{lem:VaaRH} for $\CRHDb_{\CC,0}$ gives the last assertion of Proposition \ref{prop:VCRH}.
\end{proof}

\subsubsection*{Localization and Mellin transform}
We may give similar definitions and similar arguments for the germ $\RHDbmo_{\CC,0}$. We say that a germ $\wt u\in\RHDbmo_{\CC,0}$ has order less than $\alpha'$ if it belongs to $V_{\alpha'}(\cD_{\CC,0}\wt u)$, \etc. We get in particular
\begin{eqnarray}
V'_{\alpha'}(\RHDbmo_{\CC,0})&=&\sum_{-1\leq \alpha<0}\sum_p t^{-[\alpha'-\alpha]}\cO_{\CC,0}\cO_{\ov \CC,0}[\ov t^{-1}]\cdot u_{\alpha,p},\nonumber\\
V''_{\alpha''}(\RHDbmo_{\CC,0})&=&\sum_{-1\leq \alpha<0}\sum_p \cO_{\CC,0}[t^{-1}]\ov t^{-[\alpha''-\alpha]}\cO_{\ov \CC,0}\cdot u_{\alpha,p},\nonumber\\
V_{\alpha',\alpha''}(\RHDbmo_{\CC,0})
&\defin&V'_{\alpha'}(\RHDbmo_{\CC,0})\cap V''_{\alpha''}(\RHDbmo_{\CC,0})\nonumber\\
&=&\sum_{-1\leq \alpha<0}\sum_p t^{-[\alpha'-\alpha]}\cO_{\CC,0}\cdot\ov t^{-[\alpha''-\alpha]}\cO_{\ov \CC,0}\cdot u_{\alpha,p}\label{eq:Vaamo}
\end{eqnarray}
Put 
\begin{eqnarray}
V_{\alpha',\alpha''}(\CRHDbmo_{\CC,0})
&\defin&\cC^\infty V_{\alpha',\alpha''}(\RHDbmo_{\CC,0})\nonumber\\
&=&\sum_{-1\leq \alpha<0}\sum_p t^{-[\alpha'-\alpha]}\ov t^{-[\alpha''-\alpha]}\cC^\infty_{\CC,0}\cdot u_{\alpha,p}.\label{eq:Vaamoc}
\end{eqnarray}
From \cite[Theorem 4]{B-M87} we obtain:

\begin{proposition}\label{prop:VMel}
Let $u\in \CRHDb_{\CC,0}$ and denote by  $\wt u$ its image in $\CRHDbmo_{\CC,0}$. Let $\alpha',\alpha''\in\CC$. Then $\wt u\in V_{\alpha',\alpha''}$ if and only if for any $k',k''\in\ZZ$ the poles of the meromorphic function $\cJ_{u}^{(k',k'')}(s)$ are $\leq \min(\alpha'-k',\alpha''-k'')$.\hfill\qed
\end{proposition}

Remark that it is enough to verify the previous criterion for $k''=0$ and $k'\in\ZZ$ for instance.

Put $V'_{\alpha'}(\CRHDbmo_{\CC,0})=\cup_{\beta''}V_{\alpha',\beta''}$ and define $V''_{\alpha''}$ similarly.

\begin{corollaire}
For any $\alpha',\alpha''\in\CC$ with $\alpha'-\alpha''\in\ZZ$, we have
$$
V_{\alpha',\alpha''}(\CRHDbmo_{\CC,0})=V'_{\alpha'}(\CRHDbmo_{\CC,0})\cap V''_{\alpha''}(\CRHDbmo_{\CC,0}).\eqno\qed
$$
\end{corollaire}

It also follows from Lemma \ref{lem:VaaRH}(b) and Formula \eqref{eq:Vaamo} above that
\begin{eqnarray}\label{eq:Vimage}
V_{\alpha',\alpha''}(\RHDbmo_{\CC,0})&=&{\rm image} \lefcro V_{\alpha',\alpha''}(\RHDb_{\CC,0})\longrightarrow \Dbmo_{\CC,0}\rigcro.
\end{eqnarray}
and thus a similar result for $\CRHDbmo_{\CC,0}$.

\begin{corollaire}
For any $\alpha',\alpha''\in\CC$ we have
\begin{eqnarray*}
V_{\alpha',\alpha''}(\CRHDb_{\CC,0})&=&V'_{\alpha'}(\CRHDb_{\CC,0})\cap V''_{\alpha''}(\CRHDb_{\CC,0}).
\end{eqnarray*}
\end{corollaire}

\begin{proof}
Let $u\in V'_{\alpha'}(\CRHDb_{\CC,0})\cap V''_{\alpha''}(\CRHDb_{\CC,0})$. Then $\wt u\in V_{\alpha',\alpha''}(\CRHDbmo_{\CC,0})$, hence, by \eqref{eq:Vimage}, there exists $v\in V_{\alpha',\alpha''}(\CRHDb_{\CC,0})$ such that $u-v$ is supported at $0$, \ie.\ belongs to
\[
V'_{\alpha'}(\CC[\partial_t,\partial_{\ov t}]\cdot \delta)\cap V''_{\alpha''}(\CC[\partial_t,\partial_{\ov t}]\cdot \delta)
\]
which is easily seen equal to $V_{\alpha',\alpha''}(\CC[\partial_t,\partial_{\ov t}]\cdot \delta)$.
\end{proof}

\begin{remarque}\label{rem:uut}
For $u\in \CRHDb_{\CC,0}$, we have $\alpha'(\wt u)\leq \alpha'(u)$ with equality if $\alpha'(u)<0$, and a similar result for $\alpha''$.
\end{remarque}

\subsubsection*{The morphisms $L_\alpha$}
According to (c) and (d) in Lemma \ref{lem:VaaRH}, we may give the following definition:

\begin{definition}\label{def:lalpha}
For $-1\leq \alpha\leq 0$, the linear morphism
\[
L_\alpha:\gr_{\alpha,\alpha}^{V}\CRHDb_{\CC,0}\longrightarrow \CC
\]
is obtained by taking the coefficient of the class of $\dfrac{u_{\alpha,0}}{2i\pi}$ if $-1\leq \alpha<0$ and the coefficient of the class of $\delta=-\dfrac{1}{2i\pi}\partial_t\partial_{\ov t}\log\module{t}^2$ if $\alpha=0$.
\end{definition}

It will be convenient to denote also by $L_\alpha$ the map composed with the previous $L_\alpha$ and the projection $V_{\alpha,\alpha}\rightarrow \gr_{\alpha,\alpha}^{V}$, so that $L_\alpha(u)=0$ if $\alpha'(u)<\alpha$ or $\alpha''(u)<\alpha$.

\begin{proposition}\label{prop:tdovt}
Let $v\in\gr_{-1,0}^{V}\CRHDb_{\CC,0}$ and $w\in \gr_{0,-1}^{V}\CRHDb_{\CC,0}$. We then have
\[
L_0(-\partial_t v)=L_{-1}(\ov t v)\quad\text{and}\quad L_0(-\partial_{\ov t} w)=L_{-1}(t w).
\]
\end{proposition}

\begin{proof}
Any such $v$ can be written as $\partial_{\ov t}\lefpar \sum_{p\geq 1}v_p \dfrac{u_{-1,p}}{2i\pi}\rigpar$. Then we have $L_{-1}(\ov tv)=v_1$ since $\ov t\partial_{\ov t}\log\module{t}^2=1$, and
\begin{eqnarray*}
\partial_t v=\partial_t\partial_{\ov t}\lefpar \sum_{p\geq 1} v_p\dfrac{u_{-1,p}}{2i\pi}\rigpar
\end{eqnarray*}
so $L_0(-\partial_t v)=v_1$.
\end{proof}

\begin{proposition}\label{prop:caractVDb}
For $-1\leq \alpha<0$ and $u\in V_{\alpha,\alpha}(\CRHDb_{\CC,0})$, we have
\begin{eqnarray*}
L_\alpha(u)&=&\star\res_{s=\alpha}\cJ_{u}^{(0,0)}\quad\text{with }\star\neq0.
\end{eqnarray*}
\end{proposition}

\begin{proof}
This follows from the computation in the proof of Theorem 4 in \cite{B-M87}.
\end{proof}

\begin{remarque}
For $u\in V_{0,0}\CRHDb_{\CC,0}$, we may also compute $L_0(u)$ as the residue of the Mellin transform of the \emph{localized Fourier transform} of the germ $u$:
\begin{eqnarray*}
L_0(u)&=&\star\res_{s=-1}\cJ_{\cF_\loc u}^{(0,0)}\quad\text{with }\star\neq0.
\end{eqnarray*}
\end{remarque}

\subsection{The morphism $\gr_{\alpha}^{V}C_X\cM\rightarrow C_Z\gr_{\alpha}^{V}\cM$ defined using asymptotic expansions}
Keep notation of \T\ref{subsec:vfil}. Let $\cM$ be a regular holonomic $\cD_X$-module.  In order to define morphisms $\gr_{\alpha}^{V}C_X\cM\rightarrow C_Z\gr_{\alpha}^{V}\cM$, we will show below:
\begin{assertion*}
For any open set $\Omega\subset Z$, any disc $D\subset\CC$ centered at $0$ and for $-1\leq \alpha\leq 0$, the mapping
\begin{eqnarray}\label{eq:sesquiV}
\Gamma\lefpar \Omega\times D,V_\alpha(\cM)\rigpar \ootimes_\CC\Gamma\lefpar \Omega\times D,V_\alpha(C_X\cM)\rigpar &\longrightarrow & \Gamma\lefpar \Omega,\Db_Z\rigpar \\
(m,\mu)&\longmapsto&\lefcro \varphi\mapsto L_\alpha(\langle\mu(m),\varphi\rangle)\rigcro\nonumber
\end{eqnarray}
is well defined (that it takes values in $\Db_Z$ can be seen as in \cite[lemme~1]{Barlet84}) and induces $0$ on $V_\alpha(\cM)\otimes_{\CC}V_{<\alpha}(C_X\cM)$ and $V_{<\alpha}(\cM)\otimes_\CC V_\alpha(C_X\cM)$.
\end{assertion*}
Therefore, \eqref{eq:sesquiV} well defines a $\cD_{Z,\ov Z}$-linear map
\begin{eqnarray}
\Gamma\lefpar \Omega,\gr_{\alpha}^{V}\cM\rigpar \ootimes_{\CC}\Gamma\lefpar \Omega,\gr_{\alpha}^{V}(C_X\cM)\rigpar &\MRE{\langle\cbbullet,\cbbullet\rangle_\alpha} &\Db_Z(\Omega).
\end{eqnarray}
Moreover, denoting by $N$ the action of $\partial_tt+\alpha$ on the left as well as the action of $\partial_{\ov t}\ov t+\alpha$ on the right, we have, according to Propositions \ref{prop:VRH}\eqref{VRH-5} and \ref{prop:VCRH},
\begin{eqnarray*}
\langle Nm,\mu\rangle_\alpha&=&\langle m,N\mu\rangle_\alpha.
\end{eqnarray*}
Therefore, $\langle\cbbullet,\cbbullet\rangle_\alpha$ defines, for $-1\leq \alpha\leq 0$, a $\cD_{\ov Z}$-linear morphism
\[
c_{X,\alpha}(\cM):\gr_{\alpha}^{V}C_X\cM\longrightarrow C_Z\gr_{\alpha}^{V}\cM
\]
changing $N$ into $C_Z(N)$, \ie.\ the following diagram commutes:
\begin{equation}\label{eq:mor}
\begin{array}{c}
\dpl\xymatrix{
\gr_{\alpha}^{V}C_X\cM\ar[rr]^-{c_{X,\alpha}(\cM)}\ar[d]_-{N}&&C_Z\gr_{\alpha}^{V}\cM\ar[d]^-{C_Z(N)}\\
\gr_{\alpha}^{V}C_X\cM\ar[rr]^-{c_{X,\alpha}(\cM)}&&C_Z\gr_{\alpha}^{V}\cM
}
\end{array}
\end{equation}
More precisely, $c_{X,\alpha}$ is a functorial morphism between the functors $\gr_{\alpha}^{V}C_X$ and $C_Z\gr_{\alpha}^{V}$. We may also consider $c_{\ov X,\alpha}$, defined in a similar way. We then have
\begin{eqnarray}\label{eq:ca}
c_{X,\alpha}&=&C_Z\circ c_{\ov X,\alpha}\circ C_X.
\end{eqnarray}

According to Proposition \ref{prop:tdovt}, the following diagram and its Hermitian dual analogue commute:
\begin{equation}\label{eq:mor01}
\begin{array}{c}
\dpl\xymatrix{
\gr_{-1}^{V}C_X\cM\ar[rr]^-{c_{X,-1}(\cM)}\ar[d]_-{-\partial_{\ov t}}&&C_Z\gr_{-1}^{V}\cM\ar[d]^-{C_Z(t)}\\
\gr_{0}^{V}C_X\cM\ar[rr]^-{c_{X,0}(\cM)}&&C_Z\gr_{0}^{V}\cM
}
\end{array}
\end{equation}

\begin{theoreme}\label{th:cxa}
Let $\cM$ be a regular holonomic $\cD_X$-module. Then the morphism $c_{X,\alpha}$ ($-1\leq \alpha< 0$) coincides with $c_{X,\lambda}^{\psi}$ ($\lambda=\exp(2i\pi\alpha)$) and $c_{X,0}$ coincides with $c_{X,1}^{\phi}$.
\end{theoreme}

Before proving Theorem \ref{th:cxa}, we will justify the construction of the morphisms $c_{X,\alpha}(\cM)$ by proving the assertion.

\begin{lemme}
Let $\cM$ be a regular holonomic $\cD_X$-module. Then, for any $\alpha',\alpha''\in\CC$, $k\in\ZZ$ and any sections $m\in\Gamma(\Omega\times D,V_{\alpha'}\cM)$, $\mu\in\Gamma(\Omega\times D,V_{\alpha''}C_X\cM)$ and any $\varphi\in\cD^{(n,n)}(\Omega)$, the meromorphic functions $\cJ_{u_\varphi}^{(k,0)}(s)$ have poles $\leq \min(\alpha'-k,\alpha'')$.
\end{lemme}

\begin{proof}
Denote by $K$ the support of $\varphi$ and by $p\geq 0$ the order of the distribution $\mu(m)$ on $K\times D$. For $q\in\ZZ$, the functions $(t,s)\mapsto\module{t}^{2s}t^{p+q}$ and $(t,s)\mapsto\module{t}^{2s}\ov t^{p+q}$ are $C^p$ on $\{\reel(s)>-q\}\times \CC$ and depend holomorphically on $s$. Consequently, if $\chi\in\cC_c^\infty(D)$, the function $s\mapsto\langle u_\varphi,\chi(t)\module{t}^{2s}t^{p+q}dt\wedge d\ov t\rangle$ is holomorphic for $\reel(s)>-q$.

Let $b_m$ be the Bernstein polynomial for $m$ on $K\times D$: there exists $P\in\Gamma(K\times D,V_0\cD_{X\times \CC})$ such that $b_m(\partial_tt)m=tP(x,t,\partial_x,\partial_tt)m$. By assumption, the roots of $b_m(-s)$ are $\leq \alpha'$. Denote by $b_\mu$ the Bernstein polynomial for $\mu$; the roots of $b_\mu(-s)$ are $\leq \alpha''$. Fix $k\in\ZZ$, choose $r\geq 0$ so large that $p-r<\alpha'$ and consider the polynomial
\begin{eqnarray*}
B'(\partial_tt) &=& \prod_{j=-k-r}^{-k}b_m(\partial_tt+j).
\end{eqnarray*}
We have $B'(\partial_tt)t^k\mu(m)\in t^{k+r}(V_0\cD_{\Omega\times D})\cdot\mu(m)$. Hence
\begin{eqnarray*}
B'(-s)\cJ_{u_\varphi}^{(k,0)}(s)&\equiv&\cJ_{u_\psi}^{(k+r,0)}(s)\quad \mod\cO(\CC)
\end{eqnarray*}
for some $\psi\in\cC_K^\infty(\Omega)$. But $\cJ_{u_\psi}^{(k+r,0)}(s)$ is holomorphic for $\reel(s)>p-k-r$, hence for $\reel(s)>\alpha'-k$. As the zeros of $B'(-s)$ are $\leq\alpha'-k$, we conclude that the poles of $\cJ_{u_\varphi}^{(k,0)}(s)$ are $\leq \alpha'-k$.

A similar argument for $b_\mu$ shows that the poles are $\leq \alpha''$.
\end{proof}

\begin{proof}[Proof of the assertion]
If $\alpha'$ and $\alpha''$ are $<0$, the desired assertion follows from Proposition \ref{prop:VMel} and Remark \ref{rem:uut}.

In general, one uses Proposition \ref{prop:VRH}(\ref{VRH-3}) together with Proposition \ref{prop:VCRH} to show the assertion for any pair $(\alpha',\alpha'')$.
\end{proof}

\subsection{Proof of Theorem \ref{th:cxa}}
\subsubsection*{First step}
Assume that $\cM$ is supported on $Z$. One then has
\begin{align*}
i^+\cM&=V_0\cM=\ker \lefcro t:\cM\rightarrow \cM\rigcro\\
\ov\imath^+C_X\cM&=V_0C_X\cM=\ker \lefcro \ov t:C_X\cM\rightarrow C_X\cM\rigcro
\end{align*}
On the other hand one has
\[
\ker\lefcro t:\Db_{X,Z}\rightarrow \Db_{X,Z}\rigcro=\ker\lefcro \ov t:\Db_{X,Z}\rightarrow \Db_{X,Z}\rigcro
\]
and one may identify this sheaf with $\Db_Z$ by defining, for any $\mu\in\ker t$ and $\varphi\in \cC_{c}^{n,n}(Z)$ (with $n=\dim Z$), $\la\mu,\varphi\ra=\mu(\psi dt\wedge d\ov t)$, where $\psi$ is any $C_{c}^{\infty}$ $(n,n)$ form on $X$ such that $\psi_{|Z}=\varphi$. The pairing $C_X\cM\otimes _\CC\cM\rightarrow \Db_{X,Z}$ induces a pairing
\[
V_0C_X\cM\ootimes_\CC V_0\cM\longrightarrow \ker t\simeq\Db_Z.
\]
This is the pairing constructed in Corollary \ref{cor:imdir}. It coincides with the pairing defined with the help of $L_0$ in \eqref{def:lalpha}.

The theorem being true for modules supported on $Z$, it folllows that it is enough to prove it for modules satisfying $\cM=j_+j^+\cM=\cM[t^{-1}]$. Moreover, as $(c_{X,1}^{\psi},c_{X,1}^{\phi})$ and $(c_{X,-1},c_{X,0})$ are both compatible with $\can$ and $\Var$, it is enough to prove the theorem for $-1\leq \alpha<0$.

\subsubsection*{Second step}
Assume now that $\cM=j_+j^+\cM=\cM[t^{-1}]$. We will show that the nondegenerate pairing
\[
\cH^0(\ov\imath^+\Cmz_X\cM)\ootimes_\CC\cH^0(i^+j_\dag j^+\cM)\longrightarrow \Db_Z
\]
given by Corollary \ref{cor:imdir} coincides with that defined with the help of $c_{X,0}$. Notice that $C_X\cM=\ov\jmath_\dag\ov\jmath^+C_X\cM$, so that
\begin{eqnarray*}
\cH^0(\ov\imath_+\ov\imath^+\Cmz_X\cM)&=&\ker\lefcro \loc:C_X\cM\rightarrow \Cmz_X\cM\rigcro\\
\cH^0(i_+i^+j_\dag j^+\cM)&=&\coker\lefcro \coloc:j_\dag j^+\cM\rightarrow \cM\rigcro
\end{eqnarray*}
Identify
\[
\cH^0(\ov\imath^+\Cmz_X\cM)\qqbox{with}V_0\cH^0(\ov\imath_+\ov\imath^+\Cmz_X\cM)\subset V_0C_X\cM
\]
and
\[
\cH^0(i^+j_\dag j^+\cM)\qqbox{with}V_0\cH^0(i_+i^+j_\dag j^+\cM)\quad(\text{a quotient of }V_0\cM).
\]
Let $\mu$ be a local section of $V_0C_X\cM$ and $m$ a local section of $V_0\cM$. Then, if $\mu$ is in $\cH^0(\ov\imath^+\Cmz_X\cM)$, the distribution $\mu(m)$ is supported on $Z$ and is in $\ker \ov t$. We may thus apply the first step to get the result.

\subsubsection*{Third step: proof for $-1\leq \alpha<0$}
We may assume that $\cM=j_+j^+\cM$. We will show that $c_{X,\alpha}$ can be computed using $\cM_{\alpha,p}$, as we did for $c_{X,\lambda}^{\psi}$. The second step will then give $c_{X,\alpha}=c_{X,\lambda}^{\psi}$.

We have isomorphisms
\begin{multline*}
\gr_{\alpha}^{V}C_X(\cM)\arrowr{\eqref{eq:a0p}}{\sim}\ker\ov
t\partial_{\ov t} \lefpar\subset\gr_{-1}^{V}C_X(\cM)_{\alpha,p}\rigpar\\[8pt] 
\shoveright{\arrowr{\text{Lemma \ref{lem:herm}}}{\sim}\ker\ov t\partial_{\ov t}\lefpar \subset\gr_{-1}^{V}\Cmz_X(\cM_{\alpha,p})\rigpar
\MRE{\sim}\ker\ov t\partial_{\ov t}\lefpar \subset\gr_{-1}^{V}C_X(\cM_{\alpha,p})\rigpar}\\[8pt] \arrowr{c_{X,-1}(\cM_{\alpha,p})}{\sim}C_Z(\coker t\partial_t) \lefpar \subset C_Z(\gr_{-1}^{V}\cM_{\alpha,p})\rigpar \arrowr{\eqref{eq:bp0}}{\sim}C_Z(\gr_{\alpha}^{V}\cM).\\
\end{multline*}
We will identify the composed isomorphism
\begin{eqnarray}\label{eq:isocc}
\gr_{\alpha}^{V}C_X(\cM)&\isom& C_Z(\gr_{\alpha}^{V}\cM)
\end{eqnarray}
with $c_{X,\alpha}(\cM)$.

Let $\mu$ be a local section of $V_\alpha C_X\cM$ and put
\[
\wt\mu_{\alpha,p}=\sum_{k=0}^{p}[-(\partial_{\ov t}\ov t+\alpha)]^k\mu\otimes \ov e_{\alpha,k}\in C_X(\cM)_{\alpha,p}.
\]
According to Lemma \ref{lem:herm}, we may view $\wt\mu_{\alpha,p}$ as a local section of $\Cmz_X(\cM_{\alpha,p})$ by putting, for $\sum_{\ell=0}^{p}m_\ell\otimes e_{\alpha,\ell}\in\cM_{\alpha,p}$,
\[
\Bigl\langle \wt\mu_{\alpha,p},\sum_{\ell=0}^{p}m_\ell\otimes e_{\alpha,\ell}\Bigr\rangle= \sum_{k,\ell=0}^{p}[-(\partial_{\ov t}\ov t+\alpha)]^k\mu(m_\ell)\cdot u_{-\alpha-2,k+\ell-p}.
\]
To understand the image of (the class of) $\wt\mu_{\alpha,p}$ by the morphism $c_{X,-1}(\cM_{\alpha,p})$, we fix a local form $\varphi$ of maximal degree and with compact support on $Z$ and consider, under the condition that all $m_\ell$ are in $V_\alpha\cM$, the coefficient of $\dfrac{u_{-1,0}}{2i\pi}$ in
\begin{equation}\label{eq:numap}
\sum_{k,\ell=0}^{p}\la [-(\partial_{\ov t}\ov t+\alpha)]^k\mu(m_\ell),\varphi\ra\cdot u_{-\alpha-2,k+\ell-p}.
\end{equation}
The only terms contributing to it are those for which $k+\ell=p$. Put
\[
\la\mu(m_\ell),\varphi\ra=\sum_{j\geq 0}v_{\ell,j}\dfrac{u_{\alpha,j}}{2i\pi}\qbox{with }v_{\ell,j}\in\CC.
\]
The coefficient of $u_{-1,0}$ in \eqref{eq:numap} is $\sum_{k=0}^{p}(-1)^kv_{p-k,k}$.

On the other hand we have
\begin{eqnarray*}
L_\alpha\biggl( \biggl\langle\mu \Bigl(\sum_{\ell=0}^{p}[-(\partial_tt+\alpha)]^{\ell}m_{p-\ell}\Bigr), \varphi\biggr\rangle\biggr) &=& L_\alpha\Bigl(\sum_{\ell=0}^{p}(-1)^\ell\sum_{j\geq 0} v_{p-\ell,j}\dfrac{u_{\alpha,j-\ell}}{2i\pi}\Bigr)\\
&=&\sum_{\ell=0}^{p}(-1)^\ell v_{p-\ell,\ell}.
\end{eqnarray*}
Consequently, \eqref{eq:isocc} coincides with $c_{X,\alpha}(\cM)$. This ends the proof of Theorem \ref{th:cxa}.\hfill\qed

\subsection{Relation with some results of D.~Barlet}
We will show that Theorems \ref{th:compatible} and \ref{th:cxa} give generalization to regular holonomic modules of some results of D.~Barlet concerning effective contribution of monodromy to poles of $\int\module{f}^{2s}$ for a holomorphic function $f:Z\rightarrow \CC$ on a smooth manifold $Z$ (\cf.\ \cite{Barlet84}). Remark that the assumption on monodromy made by D.~Barlet concerns monodromy on the cohomology of the Milnor fibre of $f$; here however, the assumption concerns monodromy on the complex of nearby or vanishing cycles and may give better results (see \eg.\ \cite{Malgrange83}).

Let $\psi_f\CC_Z$ and $\phi_f\CC_Z$ denote the complexes of nearby and vanishing cycles (see \cite{Deligne73}) and, for $\lambda\in\CC^*$, denote by $\psi_{f,\lambda}\CC_Z$ and (for $\lambda=1$) $\phi_{f,1}\CC_Z$ the complexes corresponding to the eigenvalue $\lambda$ (and $1$) of the monodromy (see the construction in \cite{Brylinski86} or \cite{MSaito86}). These complexes are perverse up to a shift and are equipped with a nilpotent endomorphism (the nilpotent part of monodromy). Let $M_\bbullet$ denote the monodromy filtration in the perverse category (see \eg.\ \cite[\T1.6]{Deligne80} or \cite[\T1.3.9]{MSaito86}) .

\begin{corollaire}\label{cor:barlet}
Let $x^o\in f^{-1}(0)$ and assume that $x^o$ belongs to the support of $\gr_\ell^M\psi_{f,\lambda}\CC_Z$ for some $\lambda=\exp(2i\pi\alpha)\in\CC^*$ ($-1\leq \alpha<0$) and $\ell\in\NN$. Then for any sufficiently small neighbourhood $V$ of $x^o$ there exists $\varphi\in\cD^{(n,n)}(V)$ ($n=\dim Z$) such that the function
\[
I_\varphi:s\longmapsto\int_Z\module{f}^{2s}\varphi
\]
has a pole of order at least $\ell$ at some $\alpha-k$ with $k\in\NN$. Similarly, if $x^o$ belongs to $\supp\gr_\ell^M\phi_{f,1}\CC_Z$, then for each $V$ there exists $\varphi$ such that the pole order is at least $\ell+1$ at some negative integer.
\end{corollaire}

\begin{Remarques*}
\begin{enumerate}
\item
If $\ell_0$ is the maximal integer $\ell$ such that $x^o$ belongs to $\supp\gr_{\ell}^{M}\psi_{f,\lambda}\CC_Z$ (or belongs to $\supp\gr_{\ell}^{M}\psi_{f,1}\CC_Z\cup\supp\gr_{\ell}^{M}\phi_{f,1}\CC_Z$ if $\lambda=1$), then the pole of any function $I_\varphi(s)$ at points $\alpha-k$, for $\varphi$ supported in a small neighbourhood of $x^o$, has order $\leq \ell_0$.
\item
If $I_\varphi$ has a pole of order $\ell$ at some $\alpha-p$ for some $\varphi$, then for any $p'\geq p$ there exists $\psi$ such that $I_\psi$ has a pole of order $\ell$ at $\alpha-p'$: put $\psi=\module{f}^{2(p'-p)}\varphi$.
\end{enumerate}
\end{Remarques*}

\begin{proof}[Proof of Corollary \ref{cor:barlet}]
Denote by $i_f:Z\hookrightarrow X=Z\times \CC$ the graph inclusion. Put $\cM=i_{f+}\cO_Z$. As $C_Z\cO_Z=\ov\cO_Z$ (Dolbeault lemma) we have $C_X\cM=\ov\cM$. We then get a sesquilinear pairing
\[
S:\cM\otimes_\CC\ov \cM\longrightarrow \Db_X.
\]
Let us consider first the case of nearby cycles ($-1\leq \alpha<0$). By assumption, and using Riemann-Hilbert correspondence for nearby cycles, there is a local section $m$ of $\cM$ such that $m\in V_\alpha(\cM)$, the class of $m$ in $\gr_{\alpha}^{V}\cM$ belongs to $M_\ell\gr_{\alpha}^{V}\cM$ and its class in $\gr_{\ell}^{M}\gr_{\alpha}^{M}\cM$ is nonzero at $x^o$.

As the pairing \eqref{eq:conjprim} is nondegenerate, there exists $\mu$ in $V_\alpha\cM$ such that $S([m],N^\ell[\ov\mu])\neq0$ in $\Db_Z$. This means that there exists $\psi\in\cD^{(n,n)}(Z)$ such that, if we put
\[
m=\sum_{i\geq 0}m_i\partial _{t}^{i}\delta(t-f),\quad \mu=\sum_{j\geq 0}\mu_j\partial _{t}^{j}\delta(t-f),
\]
where $m_i,\mu_j$ are holomorphic in a neighbourhood of $x^o$, the germ
\[
\sum_{i,j}\partial_{t}^{i}\partial_{\ov t}^{j}(\partial_{\ov t}\ov t+\alpha)^\ell \lefcro \int_{f=t}m_i \ov{\mu_j}\psi\rigcro
\]
in $\CRHDb_{\CC,0}$ has a nonzero coefficient on $u_{\alpha,0}$. Hence, there exist $i$  and $j$
 such that $\partial_{t}^{i}\partial_{\ov t}^{j}\int_{f=t}m_i \ov{\mu_j}\psi$ has a nonzero coefficient on $u_{\alpha,\ell}$. The result follows from the computation of Mellin transform (\cite[Theorem 4]{B-M87}).

\medskip
The assertion for $\phi$ follows from
\begin{eqnarray}
\gr_\ell^M(\phi_1\cM)_{x^o}\neq0&\Implique&\gr_{\ell+1}^M(\psi_1\cM)_{x^o}\neq0,
\end{eqnarray}
for which we briefly recall the proof. As $\cO_Z$ is a simple $\cD_Z$-module, $\cM$ is a simple $\cD_X$-module, according to Kashiwara's equivalence theorem. In particular it has neither submodule nor quotient module supported by $Z$. This implies (see \cite[Lemme 5.1.4]{MSaito86} forgetting the filtration $F$) that $\can:\psi_1\rightarrow \phi_1$ is onto and $\Var:\phi_1\rightarrow \psi_1$ is injective. From \cite[Lemme 5.1.12]{MSaito86} we deduce that for any $\ell$ we have
\[
\can\lefpar M_\ell\psi_1\rigpar \subset M_{\ell-1}\phi_1,\quad \Var \lefpar M_\ell\phi_1\rigpar \subset M_{\ell-1}\psi_1
\]
and that the induced morphisms
\[
\can:\gr_\ell^M\psi_1\rightarrow \gr_{\ell-1}^M\phi_1,\quad \Var:\gr_\ell^M\phi_1\rightarrow \gr_{\ell-1}^M\psi_1
\]
are respectively onto and injective.
\end{proof}

\begin{remarque}
In \cite{Barlet91}, D.~Barlet introduces the topological notion of ``tangling of strata'' and shows how this tangling can be detected by inspection of the order of poles of the functions $I_\varphi(s)$. This notion has the following interpretation. Assume as in \loccit.\ that for some eigenvalue $\lambda\neq1$ the support of $\psi_{f,\lambda}\CC_Z$ is a curve $\Sigma$ near $x^o$ and assume furthermore for simplicity that the germ $(\Sigma,x^o)$ is irreducible (one may easily extend what follows to the reducible case). The complex $\psi_{f,\lambda}\CC_Z$ is perverse up to a shift by $\dim Z-1$. Let $z$ be a local coordinate on the normalization of $\Sigma$. Consider the corresponding diagram of vector spaces:
\[
\xymatrix {
\oplus_\mu \psi_{z,\mu}\psi_{f,\lambda}\CC_Z\ar@/^/[r]^-{c}& \oplus_\mu\phi_{z,\mu}\psi_{f,\lambda}\CC_Z.\ar@/^/[l]^-{v}
}
\]
The left hand term corresponds to the generic fibre of the local system $\psi_{f,\lambda}$ on $\Sigma-\{x^o\}$ and $N'=v\circ c$ is the nilpotent part of the monodromy relative to $z$ of this local system. Moreover, $\coker c$ (\resp.\ $\ker v$, \resp.\ $\ker c$) is isomorphic to the generalized eigenspace with eigenvalue $\lambda$ of the cohomology of the Milnor fibre $F_{x^o}$ of $f$ at $x^o$ in maximal degree $\dim Z-1$ (\resp.\ the cohomology with compact support, \resp.\ the cohomology in degree $\dim Z-2$). As usual, $c$ and $v$ are compatible with the direct sum decomposition indexed by $\mu$ and their $\mu$-components are isomorphisms if $\mu\neq1$. Moreover, $c$ and $v$ commute with the nilpotent part $N$ of the monodromy of $f$.

The tangling phenomenon (for the eigenvalue $\lambda$) appears when the nilpotency indices of $N$ on the \emph{cohomology sheaves} of $\psi_{f,\lambda}\CC_Z$ are strictly smaller than the nilpotency index of $N$ on the \emph{complex} $\psi_{f,\lambda}\CC_Z$. The latter can be read from the pole order of functions $I_\varphi(s)$ (Corollary \ref{cor:barlet}).

This also means that the nilpotency indices of $N$ on the spaces $\psi_{z,1}\psi_{f,\lambda}\CC_Z$ and $\coker c=H^{\dim Z-1}(F_{x^o})_\lambda$ are strictly smaller than the nilpotency index of $N$ on the space $\phi_{z,1}\psi_{f,\lambda}\CC_Z$. 

This would not happen if $c$ were \emph{strict} relatively to the monodromy filtration $M(N)$. In such a case, still denoting by $M(N)$ the monodromy filtration on $\coker c$, we would have
\[
\gr_\ell^M\coker c=\coker\gr_\ell^M c
\]
and $\gr_\ell^M\phi_{z,1}\psi_{f,\lambda}\CC_Z$ would vanish as soon as $\gr_\ell^M\psi_{z,1}\psi_{f,\lambda}\CC_Z$ and $\gr_\ell^MH^{\dim Z-1}(F_{x^o})_\lambda$ do so.

More generally, as $\im c$ and $\ker v$ are stable by $N$, the tangling phenomenon would not happen if $\phi_{z,1}\psi_{f,\lambda}\CC_Z$ could be decomposed as $\im c\oplus \ker v$, which is equivalent to the property that the canonical morphism $H_c^{\dim Z-1}(F_{x^o})_\lambda\rightarrow H^{\dim Z-1}(F_{x^o})_\lambda$ (\ie.\ $\ker v \rightarrow \coker c$) is an isomorphism (or injective, or onto, as $\dim\ker v=\dim\coker c$ by duality and self-conjugation of $\psi_f\CC_Z$). When such an isomorphism occurs, there is no ``topological tangling'' in the sense of Barlet \cite{Barlet91}.
\end{remarque}

\providecommand{\bysame}{\leavevmode\hbox to3em{\hrulefill}\thinspace}

\end{document}